\documentclass[11pt]{amsart}

\usepackage[usenames]{color}
\usepackage{amsfonts,amssymb,latexsym}
\usepackage{amsthm}
\usepackage[margin=1in]{geometry}
\usepackage{graphicx}
\usepackage{amsmath}
\usepackage[font=footnotesize,labelfont=bf]{caption}
\usepackage{subfig}
\usepackage{placeins}
\usepackage{tikz}
\usepackage{marginnote}

\usetikzlibrary{calc,intersections}

\begin{document}

\numberwithin{figure}{section}
\newtheorem{thm1}{Theorem}[section]
\newtheorem{thm}[thm1]{Theorem}
\newtheorem{lem}[thm1]{Lemma}
\newtheorem*{thm2}{Theorem}
\newtheorem*{definition}{Definition}
\newtheorem{cor}[thm1]{Corollary}
\newtheorem{claim}[thm1]{Claim}
\newtheorem{conj}[thm1]{Conjecture}
\newtheorem*{rem}{Remark}
\newtheorem*{question}{Question}
\newtheorem{ex}{Example}[section]

\newtheorem{corollary}[thm1]{Corollary}
\newtheorem{conjecture}[thm1]{Conjecture}
\newtheorem{theorem}[thm1]{Theorem}
\newtheorem{lemma}[thm1]{Lemma}
\newtheorem{exer}[thm1]{Exercise}
\newtheorem{proposition}[thm1]{Proposition}
\newtheorem{prop}[thm1]{Proposition}

\newcommand{\bi}{\begin{itemize}}
\newcommand{\ei}{\end{itemize}}
\newcommand{\be}{\begin{enumerate}}
\newcommand{\ee}{\end{enumerate}}
\newcommand{\ds}{\displaystyle}
\newcommand{\ul}{\underline}
\newcommand{\rk}{{\color{red} \noindent{\bf Remark }}}
\newcommand{\bl}{\color{blue}}
\newcommand{\gr}{\color{green}}
\newcommand{\red}{\color{red}}
\newcommand\C{\mathbb{C}}
\newcommand\R{\mathbb{R}}
\newcommand\Q{\mathbb{Q}}
\newcommand\Z{\mathbb{Z}}
\newcommand\N{\mathbb{N}}
\newcommand{\hnote}[1]{\marginnote{ \scriptsize \textcolor{red}{HNH:{#1}}}}

\interfootnotelinepenalty=10000

\title{Bounding Crossing Number in Rectangular and Hexagonal Knot Mosaics}
\date{\today}
\author{Hugh Howards, Jiong Li, Xiaotian Liu}


\address{Department of Mathematics, Wake Forest University, Winston-Salem, NC 27109, USA}

\maketitle

\begin{abstract}
Howards and Kobin give a sharp upper bound for crossing number of knots on rectangular mosaics \cite{hk}.  Here we extend the  proof to create a new bound for hexagonal mosaics in all three natural settings and shorten the proof in the rectangular setting.
\end{abstract}

\section{Introduction}

Mosaics were introduced as a way of modeling quantum knots by Lomonaco and Kauffman in \cite{LK}.  Since then they have been studied in numerous papers. An $r$-mosaic in the rectangular sense of Lomonaco and Kauffman is an $r \times r$ board consisting of a combination of square tiles with arcs of a link running across them.

\begin{figure}[!h]
	\centering
	\setlength{\unitlength}{0.1\textwidth}
	\begin{picture}(10,5.5)
	\put(1.03,0.05){\includegraphics[width=.8\textwidth]{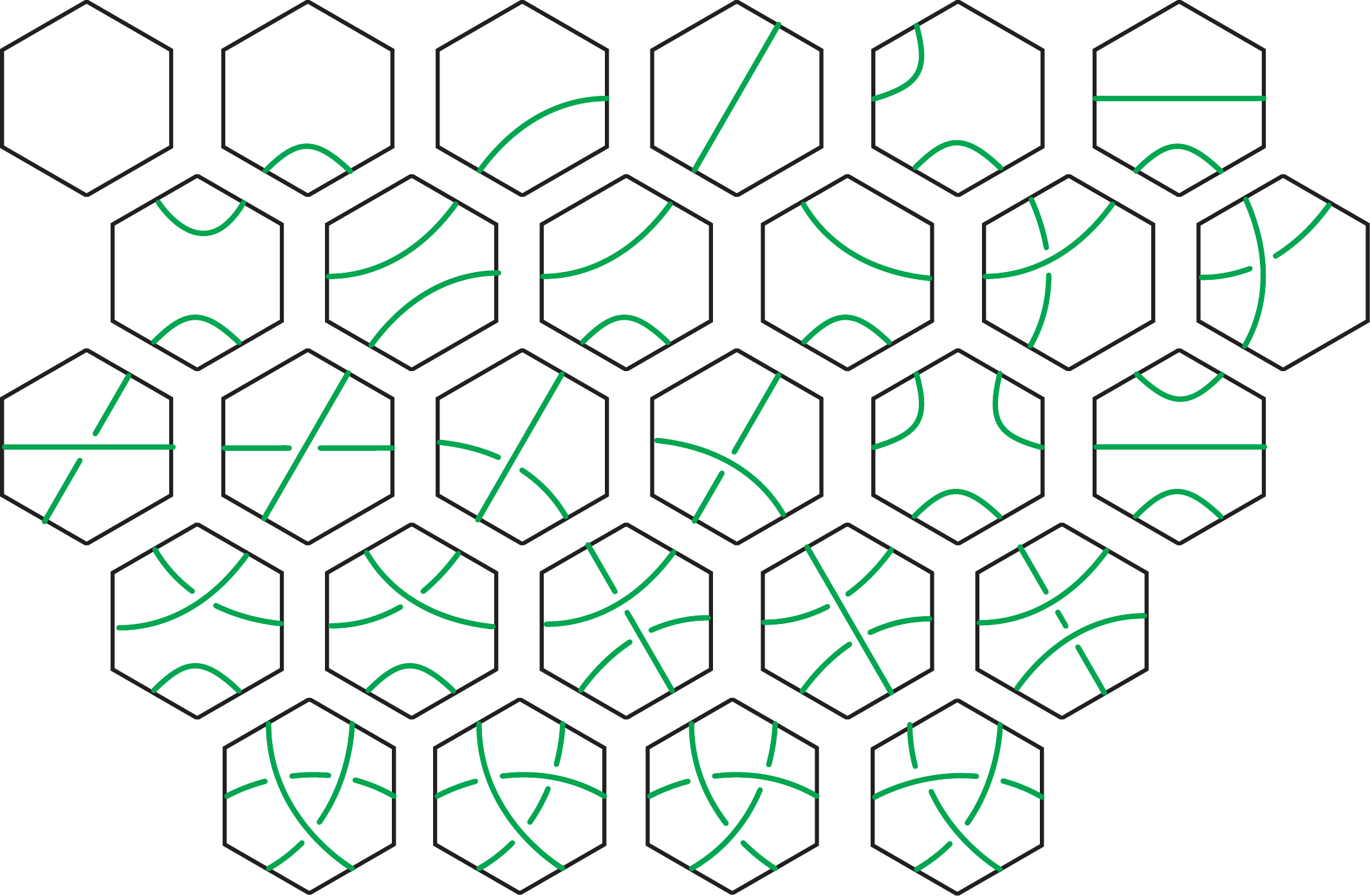}}
	\put(0.85,4.8){1}
	\put(2.11,4.8){2}
	\put(3.44,4.8){3}
	\put(4.67,4.8){4}
	\put(5.95,4.8){5}
	\put(7.2,4.8){6}
	\put(1.5,3.75){7}
	\put(2.75,3.75){8}
	\put(4.0,3.75){9}
	\put(5.2,3.75){10}
	\put(6.5,3.75){11}
	\put(7.77,3.75){12}
	\put(0.75,2.7){13}
	\put(2.09,2.7){14}
	\put(3.34,2.7){15}
	\put(4.57,2.7){16}
	\put(5.85,2.7){17}
	\put(7.13,2.7){18}
	\put(1.45,1.78){19}
	\put(2.7,1.78){20}
	\put(3.96,1.78){21}
	\put(5.22,1.78){22}
	\put(6.49,1.78){23}
	\put(2.11,0.7){24}
	\put(3.34,0.7){25}
	\put(4.57,0.7){26}
	\put(5.85,0.7){27}
	\end{picture}
	\caption{Here we see the hexagonal tiles up to rotation.  Numbers to the left of the tile are given so that we can refer to specific tiles by name when convenient.}
	\label{fig:tiles}
\end{figure}
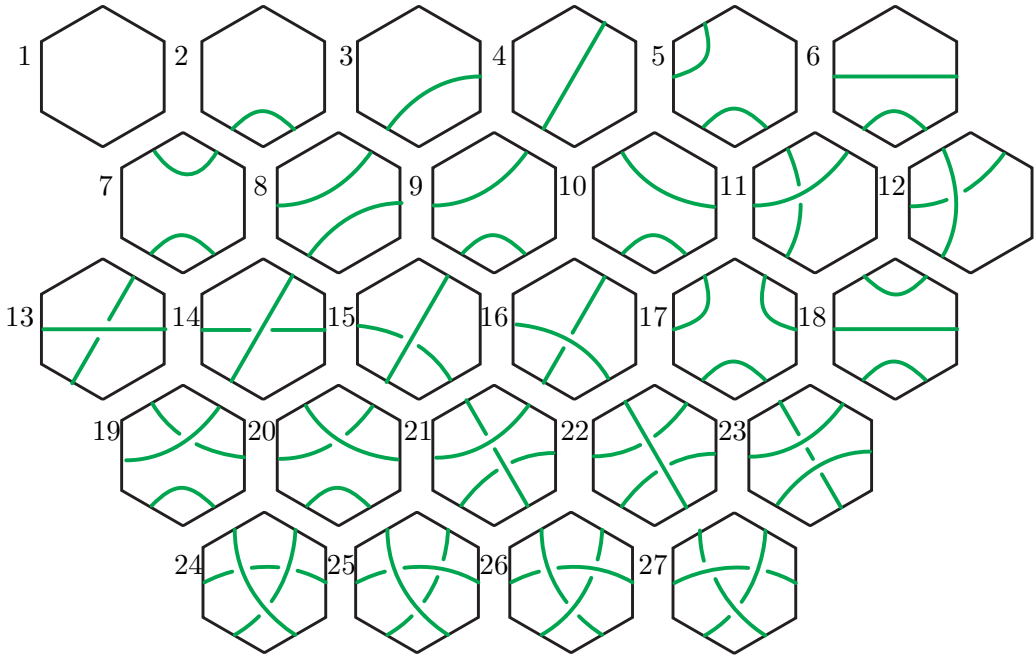

The concept of hexagonal mosaics was introduced in \cite{jmm} by McLoud-Mann et al. Instead of tiling a portion of the plane with square tiles, they tiled it with regular hexagons. The hexagonal tiles are pictured up to rotation in Figure~\ref{fig:tiles}.   
  In the hexagonal context an $r$-mosaic is the set of all tiles in a hexagonal grid within distance $r$ of a central tile.  So a hexagonal 1-mosaic contains 1 tile, a 2-mosaic contains that tile plus the 6 tiles adjacent to it and so on.  

In knot theory in general, finding bounds relating crossing number to other invariants can be notoriously difficult.  
Howards and Kobin give the following sharp bound relating crossing number and mosaic number for rectangular $r$-mosaics.  

\medskip

\noindent {\bf Theorem 8.2 \cite{hk}}. \emph{ 
	Given a rectangular $r$-mosaic with $r>3$ and any knot $K$ that is projected onto the mosaic, the crossing number $c$ of $K$ is bounded above by the following:}\\
\begin{equation*}
c \leq
\begin{cases}
(r - 2)^{2} - 2 & \quad \text{if $r = 2k + 1$}\\
(r - 2)^{2} - (r - 3) & \quad \text{if $r = 2k$.}
\end{cases}
\end{equation*}

\bigskip

The odd case is easy to prove, but the proof of the even case is fairly complicated and needed the introduction of a tool called the complement of a mosaic. 
Here we refine the complement and introduce a more streamlined proof for the theorem above.  Simultaneously we also establish new, sharp bounds for hexagonal mosaics.

We note that the complement is also useful in other settings within hexagonal mosaics including in the proof in our paper  \cite{hllfam} where we extend the results of Ludwig, Paat, etc in \cite{L} to find an infinite family of knots which cannot achieve their crossing number on the smallest hexagonal grid on which they fit.

Hexagonal mosaics are interesting in their own right, but we also point out that they generate knots very efficiently.  On a hexagonal 2-mosaic there is only one tile where crossings can occur (the other six tiles are on  the boundary and have no crossings) so the only non-trivial knots we can get on a hexagonal 2-mosaic are trefoils.  See, for example, the one shown in Figure~\ref{fig:trefoil}.  The existence of a nontrivial knot on a 2-mosaic, however, already shows some of the potential power of hexagonal mosaics since in rectangular mosaics  a 1-mosaic is too small even for an unknot, 2-mosaics and 3-mosaics only yield unknots and even up through 5-mosaics the only knots one can get are the trefoil and figure eight knot,  but hexagonal 3-mosaics already yield knots up to crossing number 19 and hexagonal 5-mosaics contain knots up through crossing number 108.  

\begin{figure}[!h]
	\centering
{\includegraphics[width=.32\textwidth]{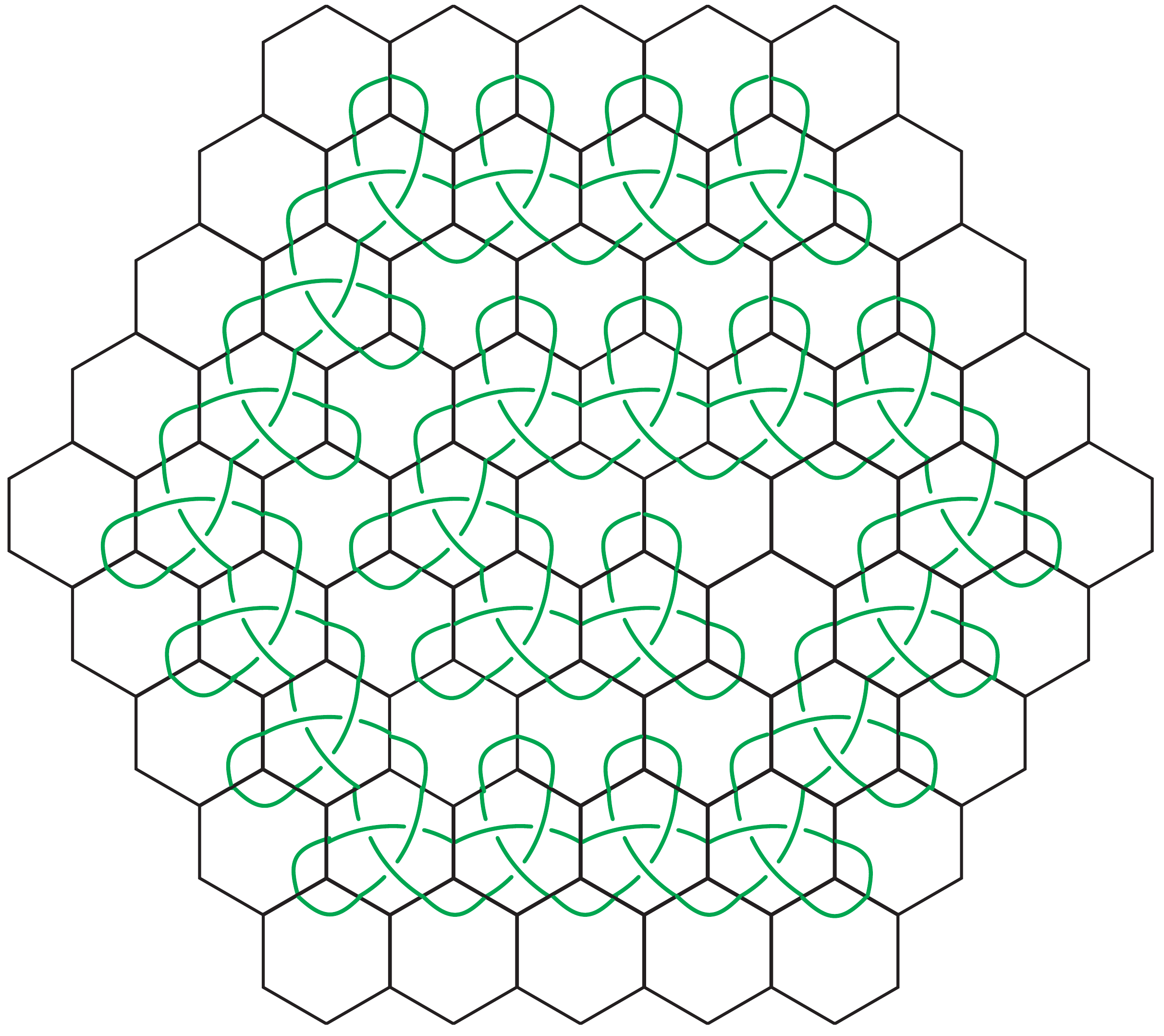}}	
{\includegraphics[width=.32\textwidth]{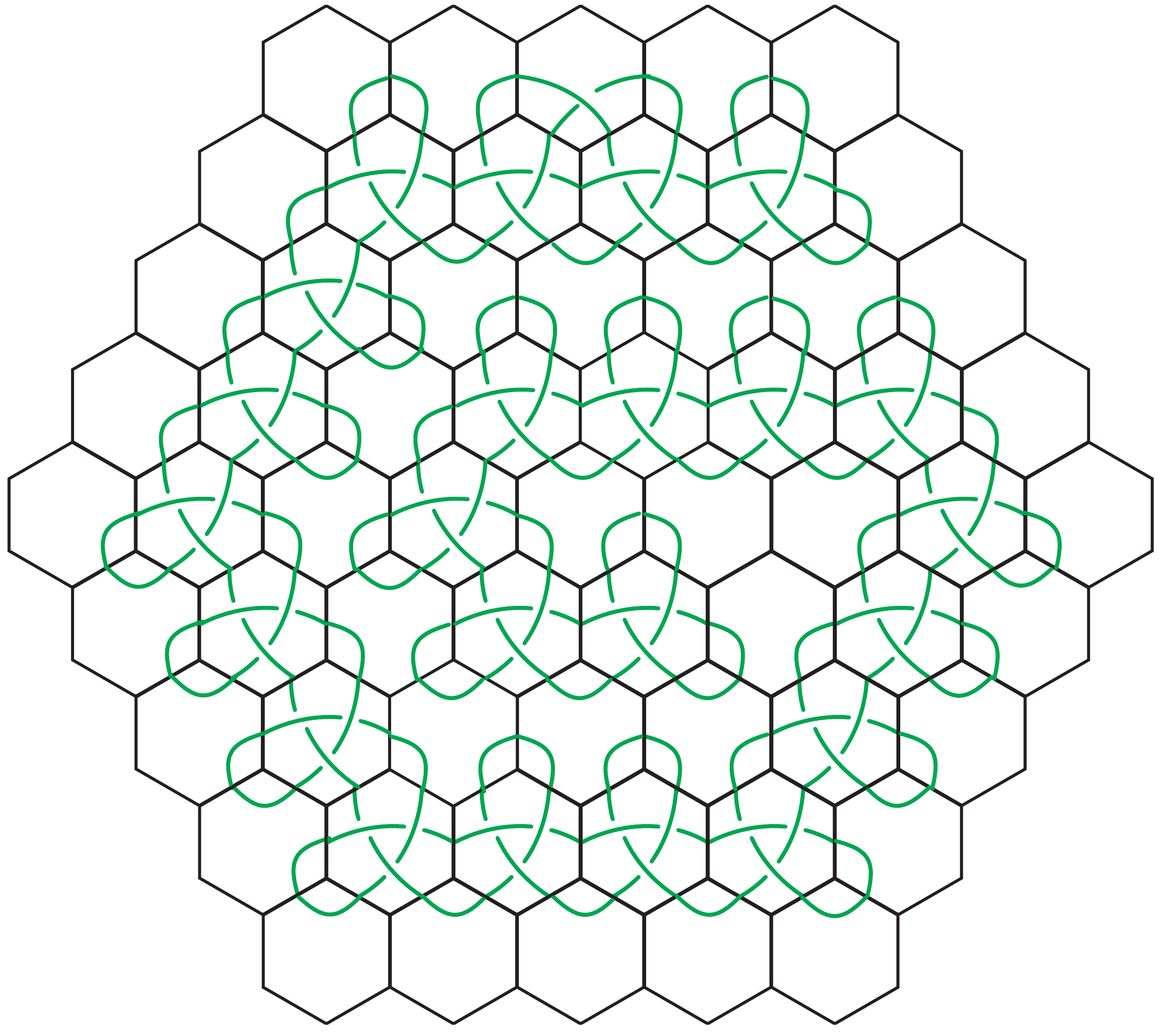}}	
{\includegraphics[width=.32\textwidth]{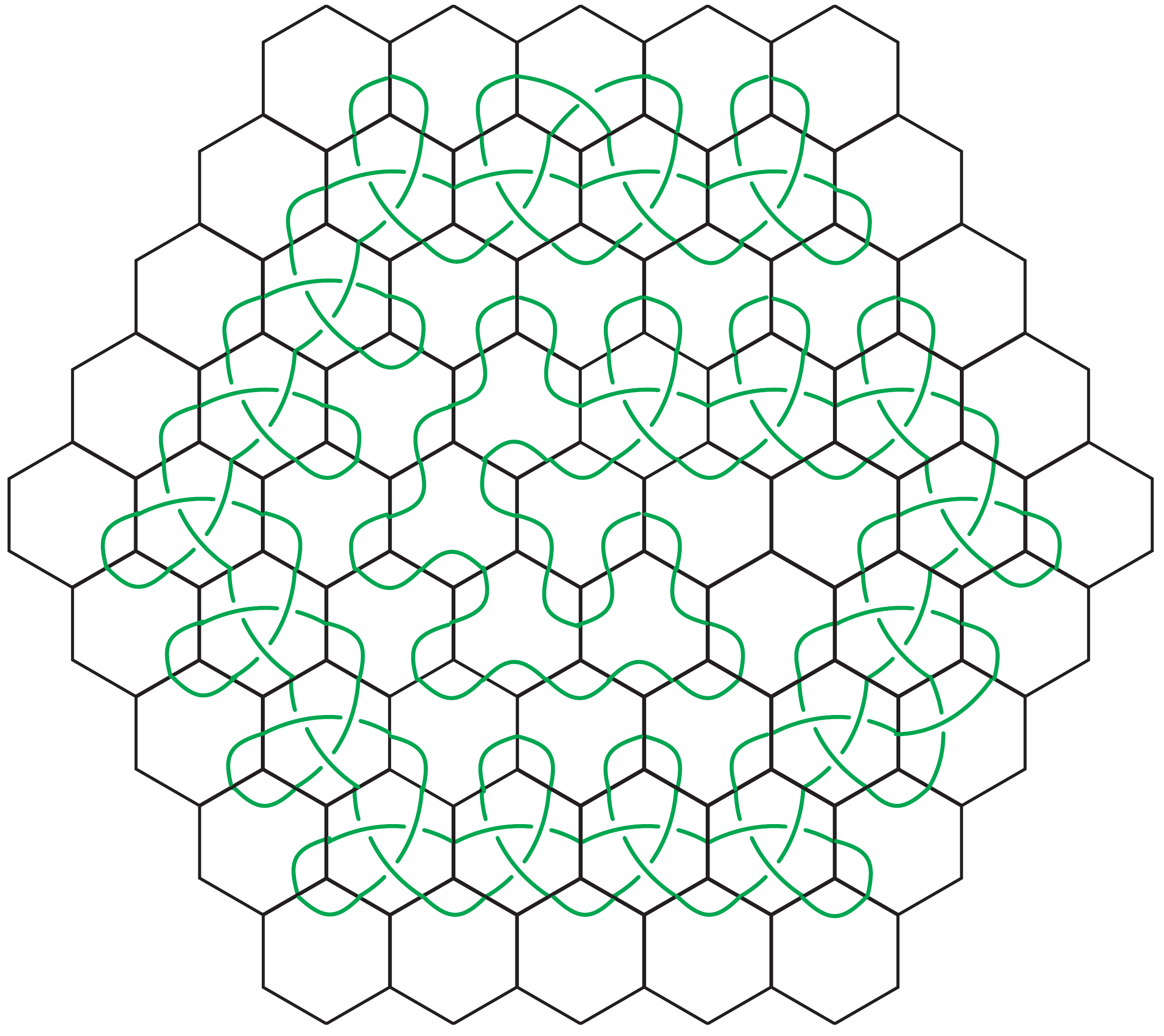}}	

	\caption{On the left we see a standard mosaic consisting of the connect sum of 23 trefoils built using only Tiles 27, 17, and 5 on the interior.  On the far right we see an enhanced mosaic, but it also shows how to lower the crossing number in any of the settings by multiples of 3 by replacing Tile 27 by Tile 17.  The center and right figures show that in the enhanced setting we can add in our choice of one or two more crossings in the boundary tiles in the enhanced setting to achieve crossing numbers of the form $3k+1$ and $3k+2$.}
	\label{fig:chain}
	\label{fig:twoways} \end{figure}

One might worry that the computational complexity is greater in hexagonal mosaics than in rectangular mosaics because 
one has more possible tiles and perhaps this cancels out the advantage of getting so many more knots on a small board, but in fact for $r=5$ where one can only get a few knots in the rectangular setting, it is easy to show that one can hit knots with every crossing number of the form $3k, k \in \Z$  for $0 \leq k \leq 23$ while using on the interior only three types of tiles: one type of crossing tile (Type 27 in Figure~\ref{fig:tiles}) together with two non-crossing tiles  (17 and 5 in the same figure), two fewer than the five options one has in the rectangular setting.  Figure~\ref{fig:chain} shows a connect sum of $k=23$ trefoils in the standard setting.  Each time we replace one of the Type 27 tiles in the figure with a Type 17 tile we drop crossing number by 3.  We do this several times in the mosaic on the right in Figure~\ref{fig:chain}.  In the enhanced setting we can add in one or two crossings in the boundary tiles allowing us to hit every crossing number of the form $3k+1$ (for $2 \leq k \leq 23$) in the middle of the figure and $3k+2$ (for $4 \leq k \leq 23$) on the right.  With a little work one can almost certainly dramatically improve these bounds and get good bounds for prime knots as well  (but we used this example because we could explain the construction for a large collection of crossing numbers without adding a large appendix of figures).  
Thus even restricting to a smaller collection of interior tiles than are available in rectangular mosaics still leads to dramatically more non-trivial knots. 

The authors appreciate the reviewer's careful reading and helpful comments. 


\section{Definitions}

\begin{figure}[!h]
	\centering
	\setlength{\unitlength}{0.1\textwidth}
	\begin{picture}(10,3.0)

	\put(0,-.057){\includegraphics[width=.32\textwidth]{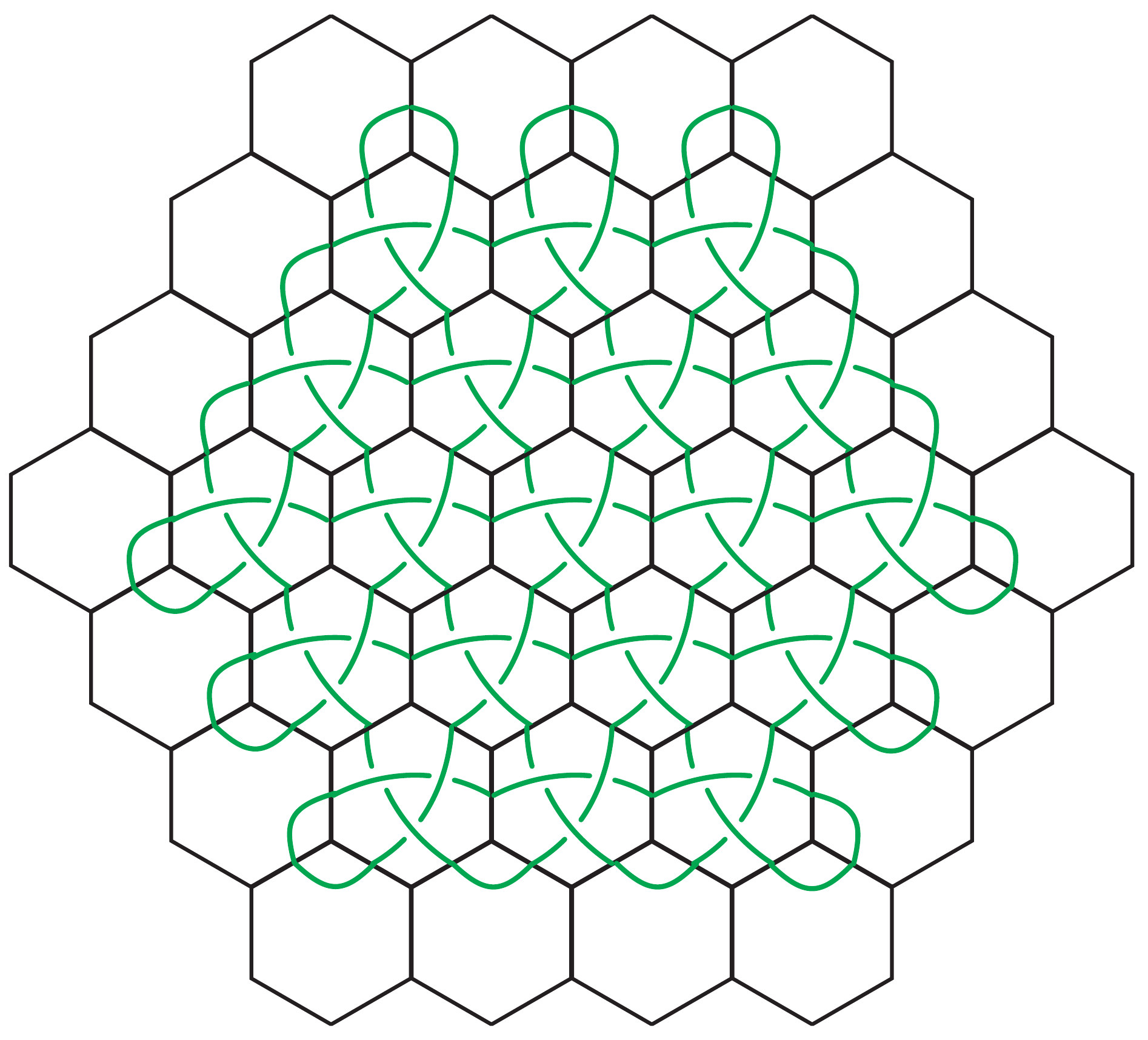}}
	\put(3.3,-.057){\includegraphics[width=.32\textwidth]{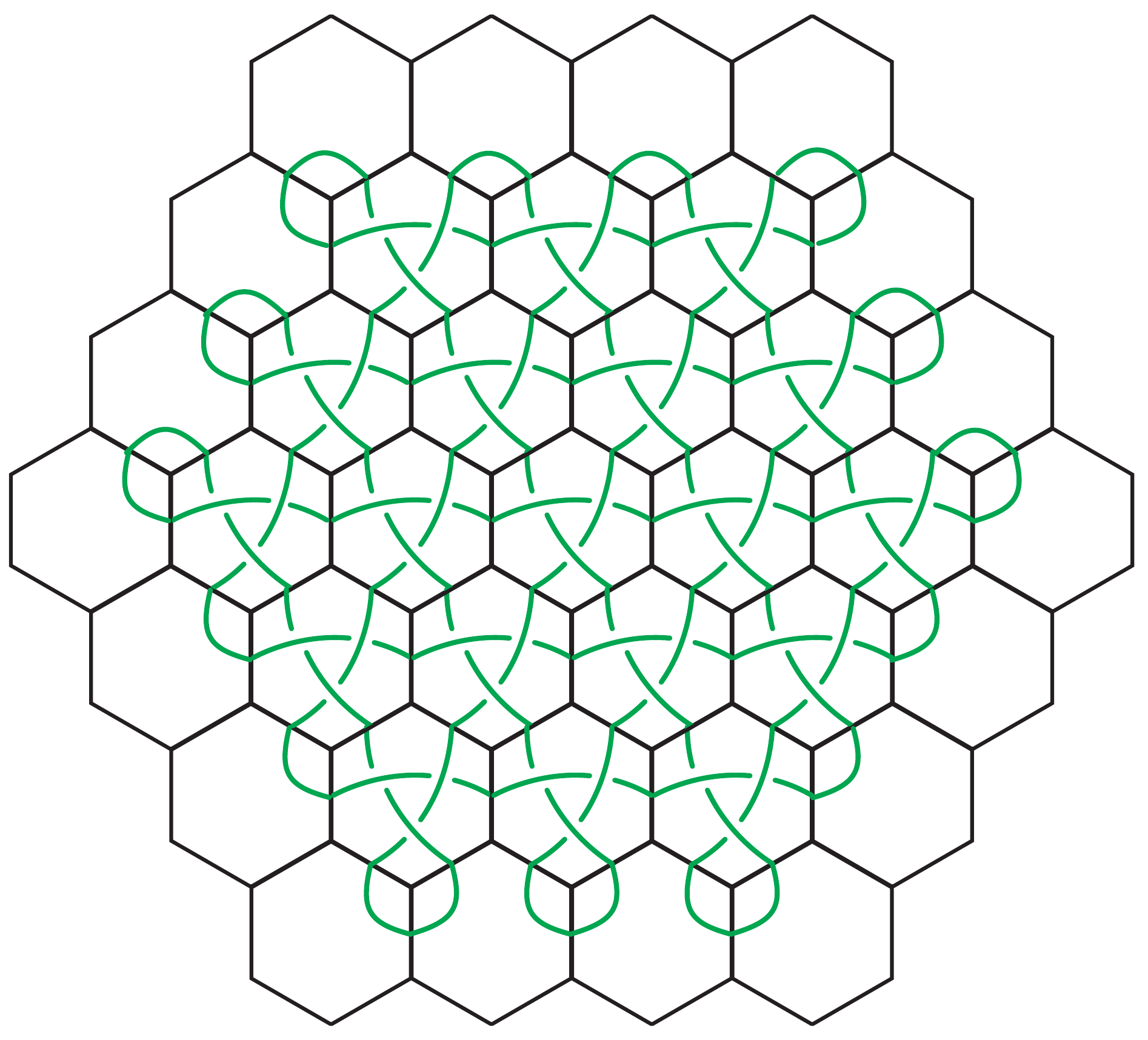}}
	\put(6.6,-.057){\includegraphics[width=.32\textwidth]{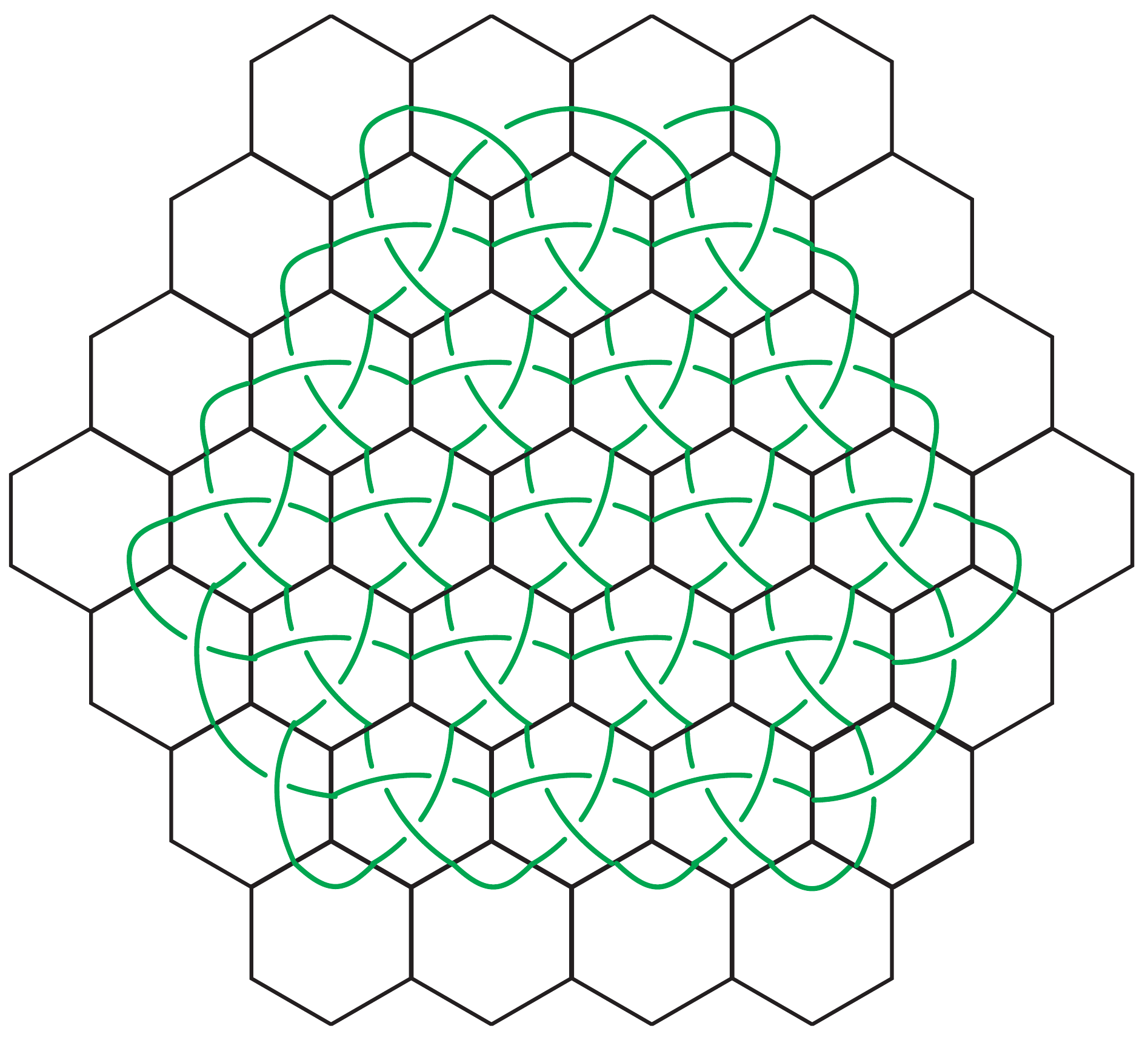}}
	\put(.3,2.5){$\ddddot{L_4}$}
	\put(6.9,2.5){$\widehat{L_4}$ }
\end{picture}
	\caption{Here we see three ways of connecting up a collection of interior tiles through the boundary.  The interior tiles on the three mosaics are identical and only the boundary tiles are different.  The link on the left is $\ddddot{L_4}$ (a standard hexagonal mosaic).  The link on the right is $\widehat{L_4}$ (an enhanced hexagonal mosaic).}
	\label{fig:l4}
	\label{fig:twoways} \end{figure}

\begin{figure}[tpb]
\centering
\includegraphics[scale=.27]{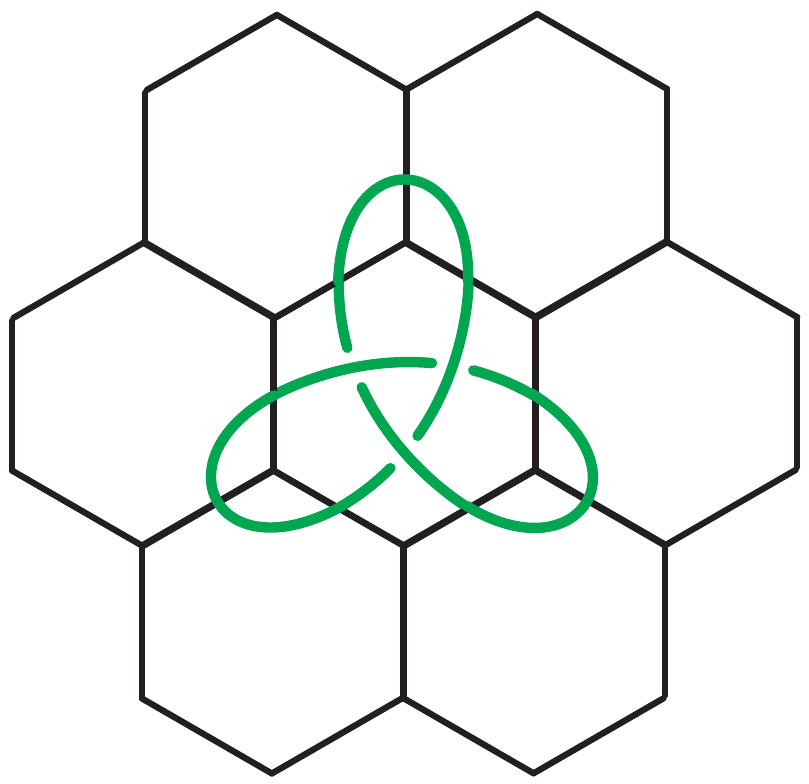}
\caption{The trefoil embedded in a hexagonal mosaic with $r=2$.  We note that for $r=2$ standard mosaics and enhanced mosaics are the same because all boundary tiles are corner tiles so it is impossible to have crossings in the boundary tiles. This single mosaic can be thought of as $\ddddot{L_2}$, $\hat{L_2}$, $\widehat{L_2}$, $\ddddot{A_2}$, $\hat{A_2}$, and $\widehat{A_2}$.}
\label{fig:trefoil}
\end{figure}

We now establish some definitions.  As mentioned above a \emph{hexagonal $r$-mosaic} is a board of radius $r$ centered around a given hexagon.    A 2-mosaic (such as the one seen in Figure~\ref{fig:trefoil}) has the original tile plus the six adjacent tiles for a total of seven tiles and so on.  We see 3-mosaics in Figure~\ref{fig:l3k3}, 4-mosaics in Figure~\ref{fig:l4}, 5-mosaics in Figure~\ref{fig:l5k5} etc. Note that a hexagonal $r$-mosaic in some sense forms a regular hexagon with $r$ tiles on each of the  outer sides just as a rectangular $r$-mosaic forms a square with $r$ tiles on each of the outer sides.
As with rectangular mosaics for a specific link on a mosaic the midpoint of each tile edge is called a {\em connection point} if it is
the endpoint of an arc drawn on that tile. Then we say that a tile within a hexagonal mosaic
is {\em suitably connected}
if each of its connection points is identified with a connection
point of an adjacent tile. As is standard practice, we only consider mosaics in which all of the tiles are suitably connected (and then we say the mosaic itself is suitably connected).

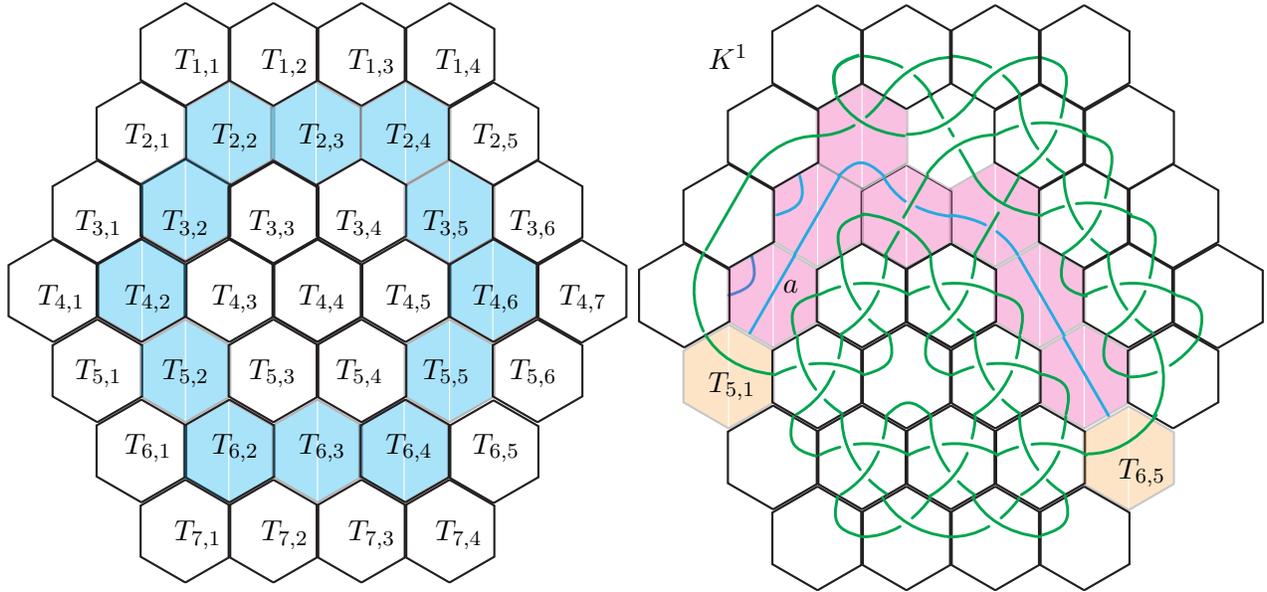
\begin{figure}[tpb]  \centering
	\setlength{\unitlength}{0.1\textwidth}
	\begin{picture}(10,5)
	\put(-0.1,.0){\includegraphics[width=.51\textwidth]{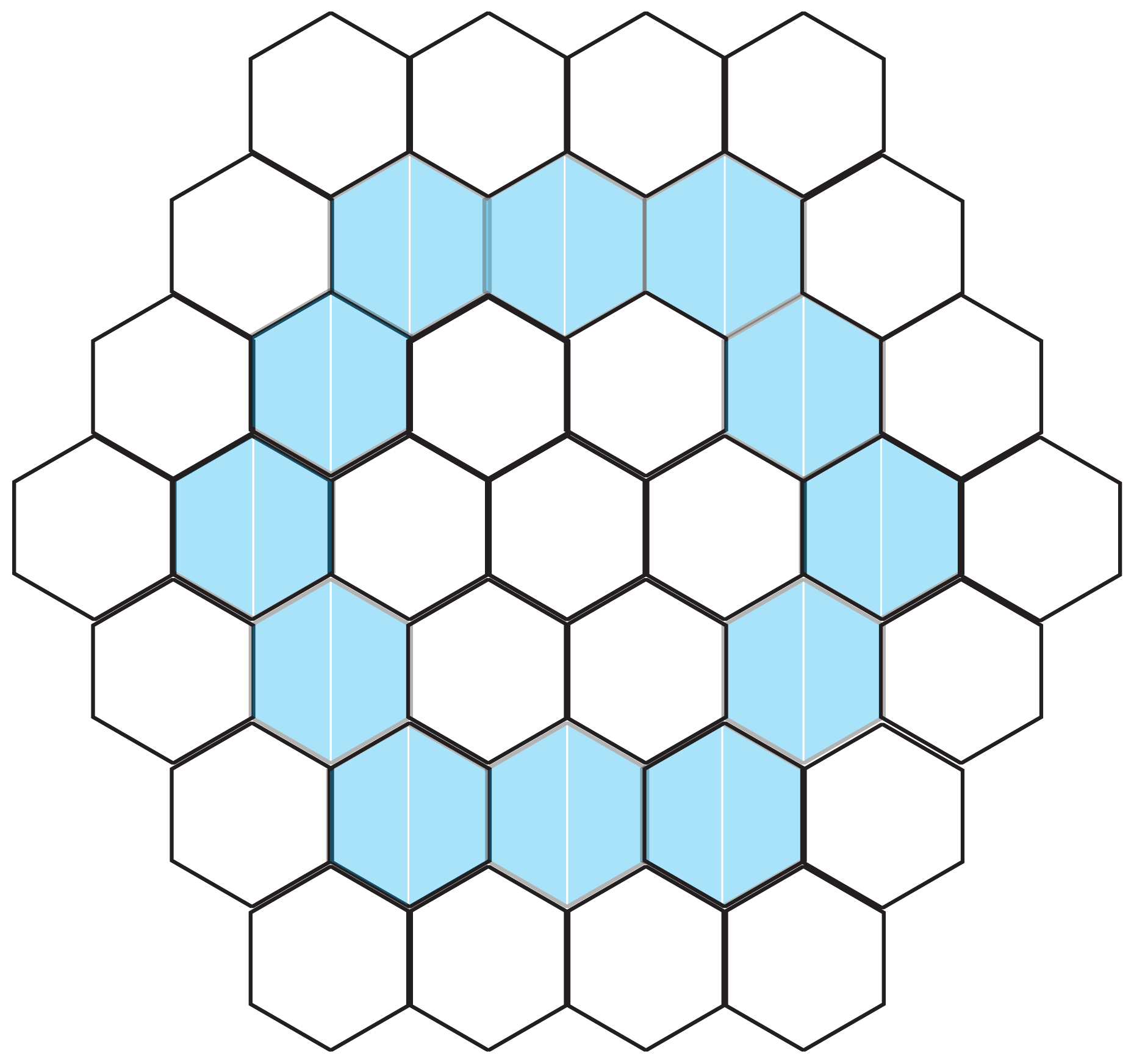}}
	\put(5,-0.06)  {\includegraphics[width=.51\textwidth]{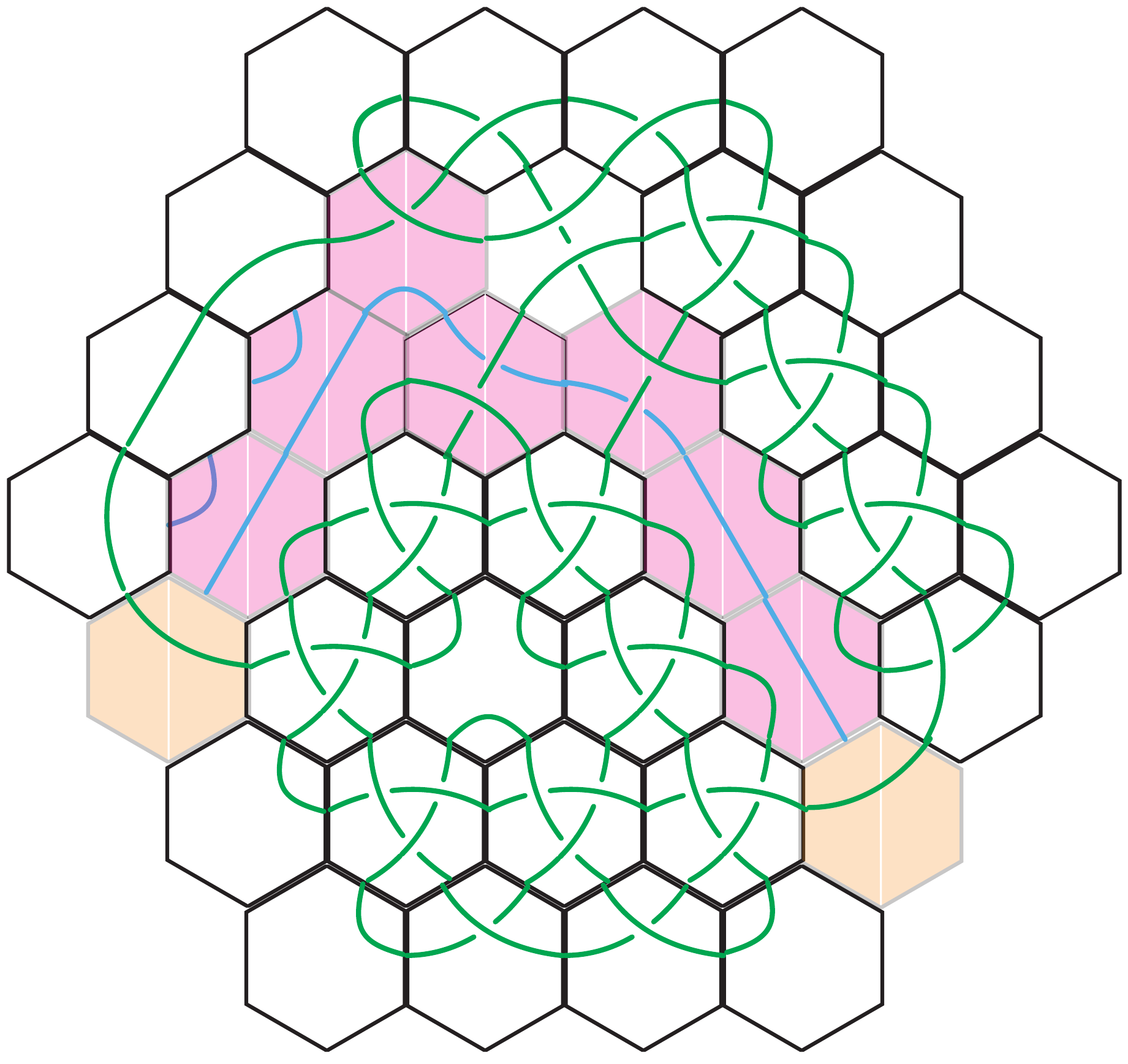}}

	\put(5.6,4.2){$K^1$}
	\put(1.3,4.2){$T_{1,1}$}
	\put(2.0,4.2){$T_{1,2}$}
	\put(2.7,4.2){$T_{1,3}$}
	\put(3.4,4.2){$T_{1,4}$}
	\put(.9,3.6){$T_{2,1}$}
	\put(1.6,3.6){$T_{2,2}$}
	\put(2.3,3.6){$T_{2,3}$}
	\put(3.0,3.6){$T_{2,4}$}
	\put(3.7,3.6){$T_{2,5}$}
	\put(0.5,2.9){$T_{3,1}$}
	\put(1.2,2.9){$T_{3,2}$}
	\put(1.9,2.9){$T_{3,3}$}
	\put(2.6,2.9){$T_{3,4}$}
	\put(3.3,2.9){$T_{3,5}$}
	\put(4.0,2.9){$T_{3,6}$}
	\put(.2,2.3){$T_{4,1}$}
	\put(.9,2.3){$T_{4,2}$}
	\put(1.6,2.3){$T_{4,3}$}
	\put(2.3,2.3){$T_{4,4}$}
	\put(3.0,2.3){$T_{4,5}$}
	\put(3.7,2.3){$T_{4,6}$}
	\put(4.4,2.3){$T_{4,7}$}
	\put(0.5,1.7){$T_{5,1}$}
	\put(1.2,1.7){$T_{5,2}$}
	\put(1.9,1.7){$T_{5,3}$}
	\put(2.6,1.7){$T_{5,4}$}
	\put(3.3,1.7){$T_{5,5}$}
	\put(4.0,1.7){$T_{5,6}$}
	\put(.9,1.1){$T_{6,1}$}
	\put(1.6,1.1){$T_{6,2}$}
	\put(2.3,1.1){$T_{6,3}$}
	\put(3.0,1.1){$T_{6,4}$}
	\put(3.7,1.1){$T_{6,5}$}
	\put(1.3,0.4){$T_{7,1}$}
	\put(2.0,0.4){$T_{7,2}$}
	\put(2.7,0.4){$T_{7,3}$}
	\put(3.4,0.4){$T_{7,4}$}
	
	\put(5.6,1.6){$T_{5,1}$}
	\put(8.9,.9){$T_{6,5}$}
	
	\put(6.2,2.4){$a$}
	
	\end{picture}
	\caption{On the left we see a system for numbering the tiles in a mosaic with the first index based on the row and the second based on how far from the left the tile lies in the row.  We also see that the outermost corona 
	consists of the boundary tiles, the penultimate corona is shaded in blue, and the central tiles are the tiles inside the shaded ring.  On the right we see a knot $K^1$ together with an arc of the complement $a$.  Our goal is to replace $K^1$ by a new knot $K^3$ which has at least as many crossings as $K^1$, but fewer arcs in its complement.}
	\label{fig:KA}
	\label{fig:grid}
	\label{fig:int}
\end{figure}

The tiles $t$ away from the center tile are called the \emph{$t^{\mathrm{th}}$ corona}, so the six tiles adjacent to the center tile on a hexagonal mosaic form the first corona and the 12 tiles just outside those form the second corona, etc.
As in the rectangular case, 
we call the tiles farthest from the center (the $r-1^{\mathrm{st}}$ corona) the \emph{boundary tiles}.    All the other tiles of the mosaic are called  \emph{interior tiles}.   We will always arrange our mosaics so that the top and bottom edges of each corona are horizontal.  So that we can specify individual tiles we will index them with the first number indicating the row from top to bottom and the second index increasing as we count tiles from left to right, so tile $T_{1,3}$ would be the third tile from the left in the top row.  See Figure~\ref{fig:KA} for a sample board with all of the tiles labeled.  
The boundary tiles (like every corona) form a hexagon of sorts.  We call the six tiles that correspond to vertices of the hexagon \emph{corner tiles} (tiles $T_{1,1}, T_{1,4}, T_{4,1}, T_{4,7}, T_{7,1},$ and $T_{7,4}$ in Figures~\ref{fig:twoways} or ~\ref{fig:grid}, for example) and the rest of the tiles that lie on the edges of the hexagon \emph{edge tiles} (such as tiles $T_{1,2}$, $T_{1,3}$ and $T_{2,5}$  in the same figures).

 At times we will also want to divide the interior tiles into their own outer ring and interior.  We call the subset of the interior tiles up through the $r-3^{\mathrm{rd}}$ corona the \emph{central tiles}
and call  the $r-2^{\mathrm{nd}}$ corona that forms the outer ring of the interior tiles the \emph{penultimate corona}.  See Figure~\ref{fig:int} for a picture demonstrating the tiles corresponding to each of the above definitions in this paragraph.   

All of the possible hexagonal mosaic tiles up to rotation are shown in Figure~\ref{fig:tiles}.   Rectangular mosaics can never have crossing tiles as boundary tiles and be suitably connected, but hexagonal mosaics in theory could.  If one prohibits crossings and only allow tiles 1 through 5 from Figure~\ref{fig:tiles} for the boundary tiles, then we call it a \emph{standard hexagonal mosaic}.  
This setting is parallel to the rectangular case in that once the interior tiles are determined there are exactly two ways the boundary tiles can connect up.   If we choose not to limit the tiles available for use as boundary tiles and thus to allow crossings in the boundary tiles we call these \emph{enhanced hexagonal mosaics}.   Finally, there is also a setting that falls between standard hexagonal mosaics and enhanced ones.  If one loosens up the boundary requirements to allow Tile 6 from Figure~\ref{fig:tiles}, but not enough to allow crossings we will call these \emph{semi-enhanced hexagonal mosaics}.
To help keep track of the setting during proofs we will denote a mosaic $M$ if it applies to all settings, $\overline{M}$ if it is on a rectangular board, $\ddddot{M}$ on a standard hexagonal board, $\hat{M}$ for semi-enhanced hexagonal,  and $\widehat{M}$ for enhanced hexagonal boards. 

A \emph{saturated} rectangular mosaic has a crossing on every interior tile.  A \emph{saturated} standard or semi-enhanced hexagonal mosaic  is one in which every interior tile is a 3 crossing tile (one of the tiles labeled 24, 25, 26, and 27 in Figure~\ref{fig:tiles}).  We will see in Lemma~\ref{lemma:boundarycrossings} that in an enhanced mosaic if the interior tiles are saturated then crossings can only occur on half of the boundary edge tiles (the boundary corona alternates between edges with crossing tiles and edges without).  Thus a \emph{saturated} enhanced hexagonal mosaic is one in which every interior tile is a 3 crossing tile and half of the edges of the boundary tiles are filled with crossings as seen in $\widehat{L_4}$ in Figure~\ref{fig:l4} and
$\widehat{L_5}$ in Figure~\ref{fig:l5satk5sat} (note: this is called a super saturated mosaic in \cite{jmm}).  
Given a link mosaic call the arcs that result from intersecting the link with the interior tiles the \emph{interior arcs} of the mosaic.

 Our results establish new bounds for standard, enhanced, and semi-enhanced hexagonal mosaics and reprove the rectangular result more concisely.   In  \cite{jmm} they require that there are no crossings on boundary tiles, but also in their proofs they require that there are only two ways to connect up through the boundary  thus they implicitly fall into the standard case and not the semi-enhanced setting.  Since the arguments would get redundant if we gave detailed proofs in all three hexagonal settings we will mainly focus this paper on standard hexagonal mosaics to match  \cite{jmm} and on the broadest case of enhanced mosaics, and will only sketch some of the details in the semi-enhanced case, but the arguments in those cases are no harder.

In \cite{hk}  the idea of the complement for a rectangular mosaic was introduced.  We extend the idea to hexagonal mosaics, but simplify it slightly as we do not need the Type 0 tiles from that paper in this new argument even in the rectangular setting.  As in \cite{hk} the complement is used to find a sharp bound on the crossing number of knots on mosaics.  In
Theorem~\ref{thm:rect} we reprove that result and in Theorem~\ref{thm:maxcross} we extend the result to hexagonal $r$-mosaics.


We now define \emph{the complement of a mosaic}.  The complement is only defined on the interior tiles of the mosaic and does not live on the boundary tiles.  Let $L$ be an embedding of a link on a  mosaic. Examine the interior tiles. On each of these tiles we pick arcs for the complement based on the arcs of $L$ already on that tile.

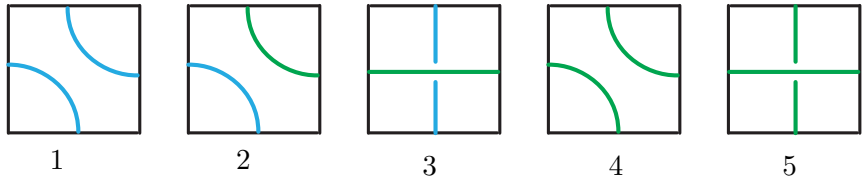
\begin{figure}[tpb]  \centering
	\setlength{\unitlength}{0.1\textwidth}
	\begin{picture}(7,1.5)
	\put(0.2,.27){\includegraphics[width=.7\textwidth]{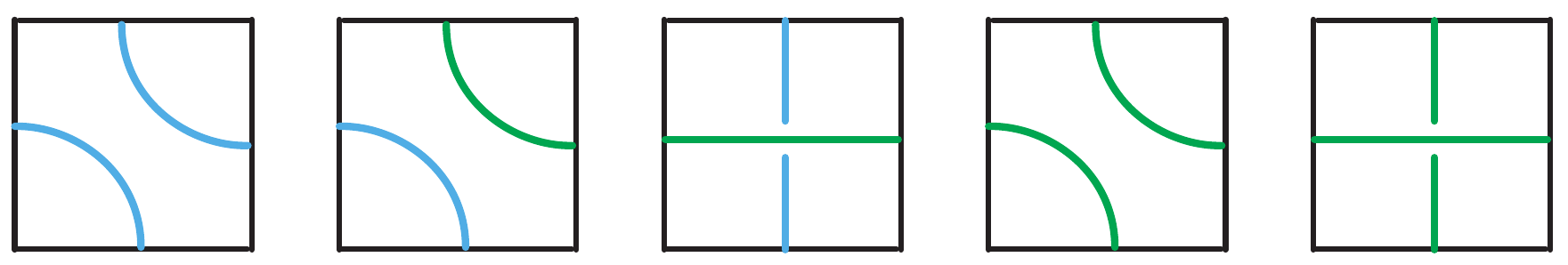}}
	\put(0.6,.05){$1$}
	\put(2.1,.05){$2$}
	\put(3.6,.0){$3$}
	\put(5.1,.0){$4$}
	\put(6.5,.0){$5$}
	\end{picture}
		\caption{The figure shows the knot intersecting a square tile
		in green arcs and the complement  is drawn in blue.  In \cite{hk} the complement for Tile 5 was blank, but Tile 4 had a single dot on it and was called a type 0 tile.  The type 0 tile is not needed in this paper and for both Tiles 4 and 5 the complement is trivial since both of those tiles contain 2 arcs of the link hitting all 4 possible connection points for the  square tile.}
	\label{fig:squarecomp}
\end{figure}

The arcs of the complement in the rectangular case are dictated according to the tiles pictured in Figure~\ref{fig:squarecomp}.   Note that if the tile is disjoint from the link as shown in Tile 1, one can pick the complement as drawn or rotated $ \frac{\pi}{2}$.  This, of course, means the complement is not well defined.  For both Tiles 4 and 5 the complement is blank (in \cite{hk} the complement for Tile 5 was blank, but the complement to Tile 4 had a dot on it to distinguish it and to indicate what was called a Type 0 tile in the complement).

\begin{figure}[tpb]
	\centering
	\includegraphics[scale=.4]{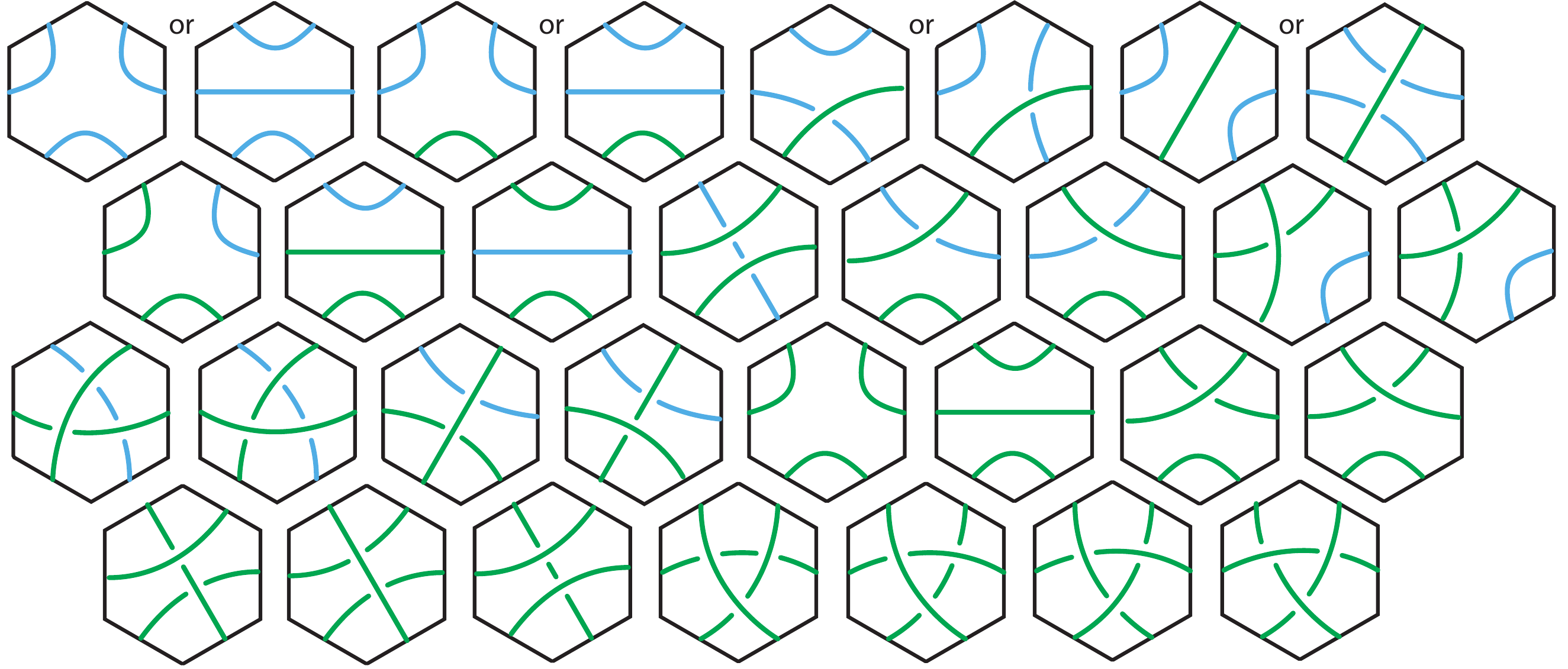}
	\caption{The figure shows the knot intersecting a hexagonal tile
		in green arcs and the complement  is drawn in blue.  When the link intersects the tile in a single arc there are usually two choices for the complement (pictured in the top row).  For the bottom row and the right half of the row above it the complement is trivial since all of those tiles contain 3 arcs of the link hitting all 6 connection points for the hexagonal tile.}
	\label{fig:comp}
\end{figure}

The choice for the arcs of the complement in the hexagonal case is dictated by the tiles pictured in Figure~\ref{fig:comp}.  Parallel to the rectangular case, the complement is defined so that it and the link together hit each the tile's potential connection points exactly once per tile.  Note that again as the figure shows the complement is not well defined.  There are four tiles where there are two different ways to choose the complement. These are the tiles where the $L$ is disjoint from the tile and the ones where it only passes through in a single arc.  Also as in the rectangular case, for the tiles where $L$ is completely disjoint from the tile, we can pick the arcs of the complement as drawn or rotated (in this case by $\frac{\pi}{3}$).   

The complement is completely contained on the interior of the mosaic and while the original link, of course, consists of a union of (possibly knotted) loops, the complement will consist of the union of $s \geq 0$ loops together with $w \geq 0$ arcs which have their end points on the boundary of the interior of the mosaic (as demonstrated by the arc of the complement drawn in Figure~\ref{fig:KA}).  As you can see in the figures defining the complement, any time an arc of the complement crosses an arc of the original link our convention dictates that the strand from the complement goes under the strand of the link. Also note that arcs of the complement never cross each other, so loops in the complement are always unknots which are split from the actual link.  

\section{Families of links denoted $L_r$ and a families of knots denoted $A_r$}


\label{sec:fam}
In this section we give the construction for a family of knots and links which will achieve the sharp bound on crossing number on $r$-mosaics.  Recall we will call a mosaic $M$, $\overline{M}$ if it is on a rectangular board, $\ddddot{M}$ on a standard hexagonal board, $\hat{M}$ on a semi-enhanced hexagonal board, and $\widehat{M}$ if it is on an enhanced hexagonal board.

\medskip
\emph{The definitions of $\overline{L_r}$ and $\overline{A_r}$ on rectangular $r$-mosaics:}
\medskip

\begin{figure}[tpb]
	\centering
	\includegraphics[scale=.25]{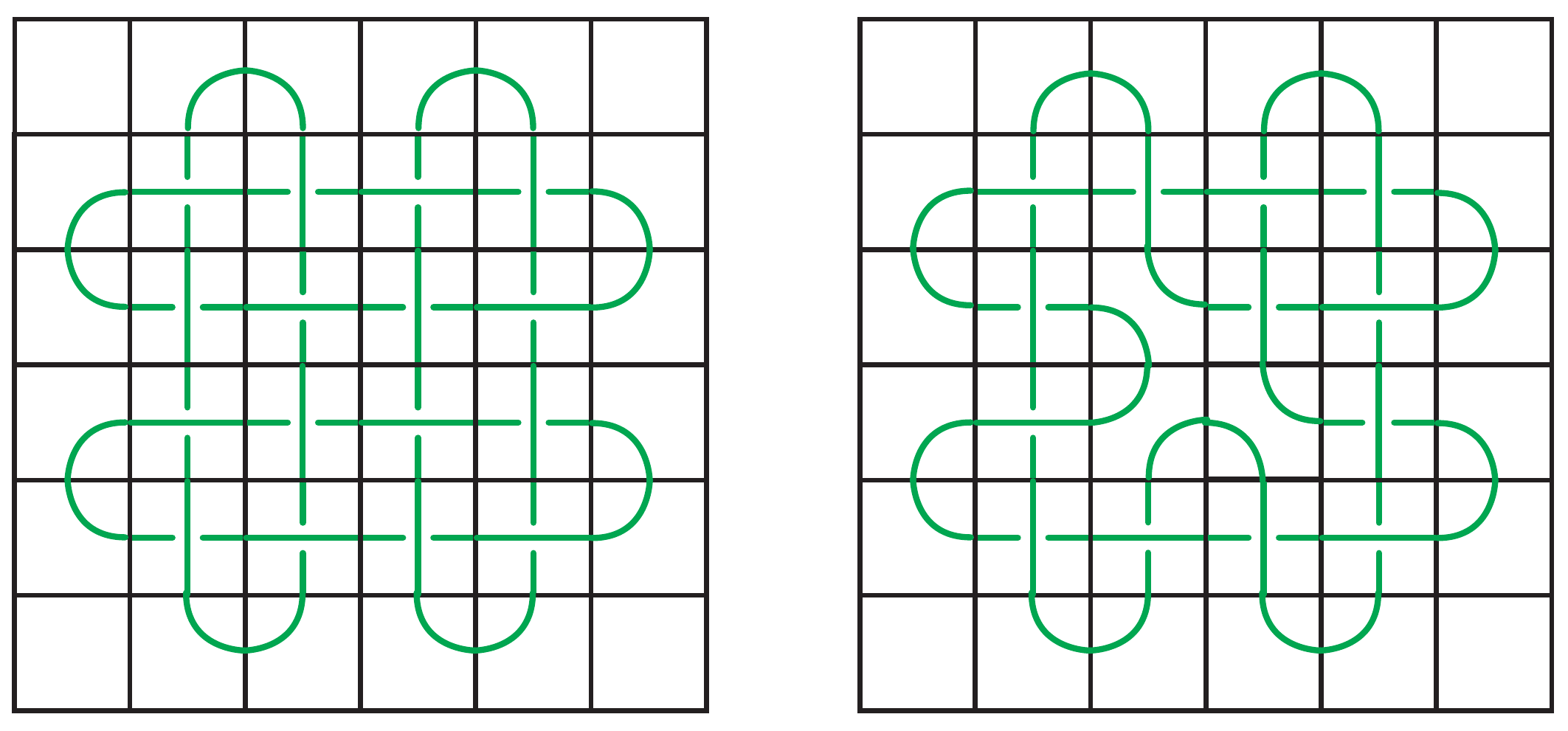}
	\caption{Here we see rectangular 6-mosaics $\overline{L_6}$, a 4 component link and $\overline{A_6}$, a maximal crossing number knot obtained by smoothing three crossings of $\overline{L_6}$.}
	\label{fig:sqfam}
\end{figure}

For rectangular mosaics the knots which achieve the sharp bound are provided in \cite{hk} and break down into distinct cases of $r$ odd or even. The theorem is easy in the odd case (say $r=2k+1, k \in \Z, k \geq 1$) so we will skip it in this paper.   
For $r$ even (let $r=2k, k \in \Z, k \geq 2$), for alternating saturated boards  once the interior tiles are chosen the two ways to choose arcs in the boundary tiles to close up the link indeed yield distinct links.  One is not reduced, having a possible type I Reidemeister move in each of the four corners.  The other is reduced and has $r-2$ components.  Call this  second, reduced link $\overline{L_r}$ (see Figure~\ref{fig:sqfam} for $\overline{L_6}$).  We can smooth $r-3$ crossings of $\overline{L_r}$ to yield a
reduced alternating knot which we now call $\overline{A_r}$.  Note that $\overline{A_r}$ is not uniquely determined.  Depending on which crossings we smooth we can get different choices for $\overline{A_r}$, but this does not cause any problems.  It is easy to pick crossings to smooth that yield a reduced knot and all of these choices for $\overline{A_r}$ on an even board have $(r-2)^2-(r-3)$ crossings and maximize crossing number. See Figure~\ref{fig:sqfam} for a picture of $\overline{A_6}$. Although not called $\overline{A_r}$ in the previous paper, these are the same knot mosaics used in \cite{hk}.

\medskip
\emph{The definitions of $\ddddot{L_r}$ and $\ddddot{A_r}$ on standard hexagonal $r$-mosaics and $\hat{L_r}$ and $\hat{A_r}$ on semi-enhanced hexagonal $r$-mosaics:}
\medskip

Now we define $\ddddot{L_r}$ and $\ddddot{A_r}$ in the standard hexagonal setting.  Unlike rectangular mosaics, the definitions are the same for both even and odd $r$ with the only exception being $\ddddot{A_3}$ as we will see below.  
Examine a  saturated hexagonal $r$-mosaic in which every interior tile is chosen to be the same alternating 3 crossing tile - either Tile 26 or 27 from Figure~\ref{fig:tiles}.  Choosing one tile yields a mirror image of the result from the other, and we will arbitrarily choose to use tile 27 from Figure~\ref{fig:tiles} in the construction in this paper.

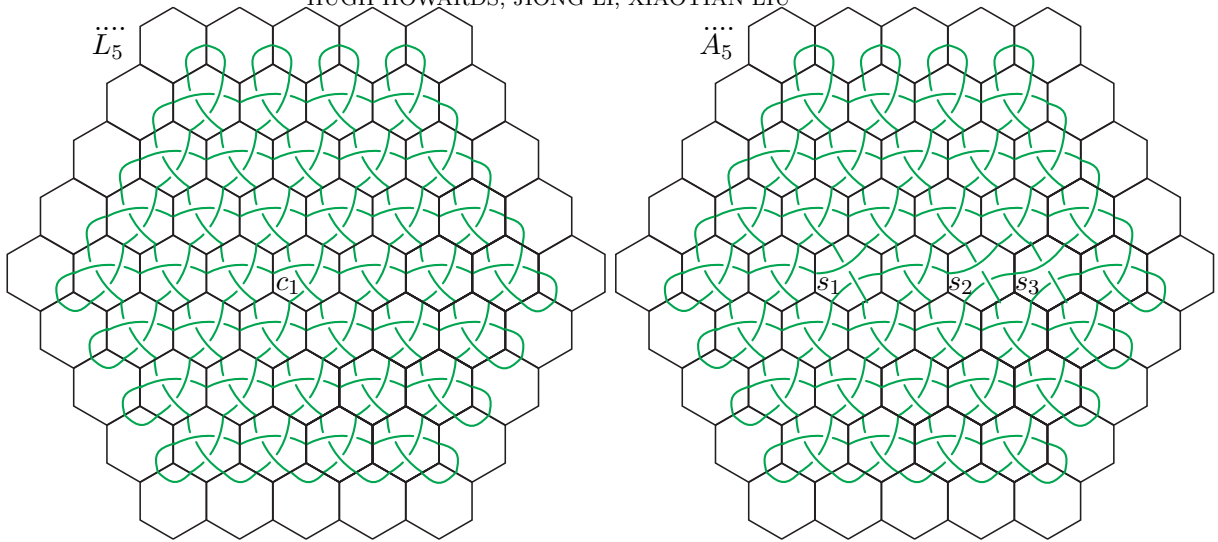
\begin{figure}[tpb]  \centering
	\setlength{\unitlength}{0.1\textwidth}
	\begin{picture}(9,4)
	\put(0.1,-.057){\includegraphics[width=.49\textwidth]{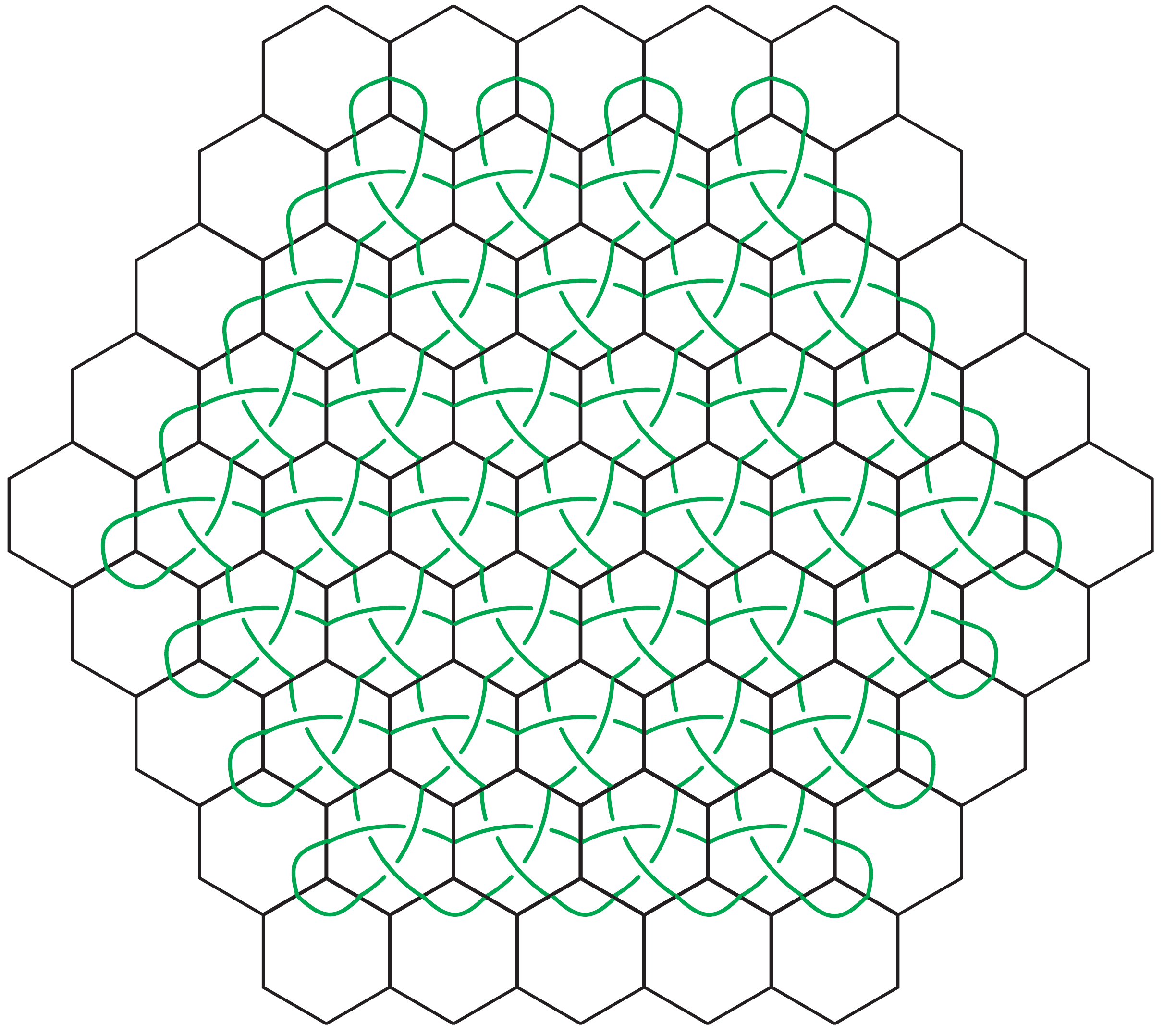}}
	\put(5,-.057){\includegraphics[width=.49\textwidth]{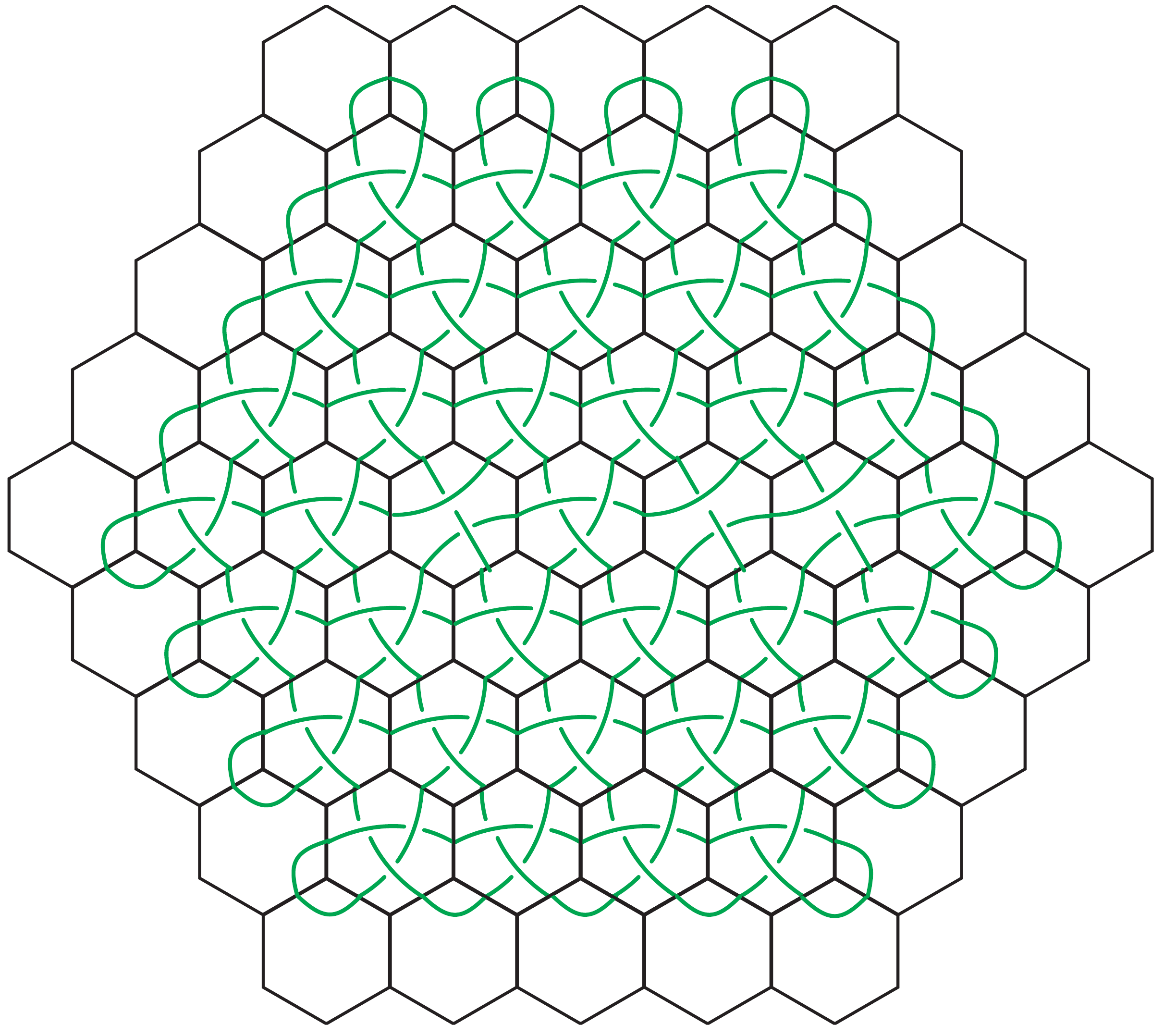}}
	\put(0.8,3.9){$\ddddot{L_5}$}
	\put(5.7,3.9){$\ddddot{A_5}$}
		\put(6.65,2.0){$s_1$}
	\put(7.71,2){$s_2$}
	\put(8.25,2){$s_3$}
	\put(2.3,2){$c_1$}

	\end{picture}
	
\caption{Here $\ddddot{L_5}$ in the standard hexagonal setting (left) is converted into an alternating knot $\ddddot{A_5}$ (right)
by smoothing crossings replacing three crossing tiles with two crossing tiles $s_1$, $s_2$, and $s_3$.  The same general strategy is used for $\ddddot{A_r}$ for all $r \geq 4$.}
\label{fig:l5k5}
\end{figure}

Once the alternating interior tiles are set the choice of the standard setting means there are only two ways to connect the link up through the boundary tiles (the two mosaics on the left in Figure~\ref{fig:l4} show the two ways one could connect up an alternating, saturated standard hexagonal 4-mosaic).  For each $r$, one of the two ways results in a reducible link with many nugatory crossings, but the other is not reducible.  The irreducible link is called $\ddddot{L_r}$.  $\ddddot{L_2}$ (which is also $\ddddot{A_2}$ since $\ddddot{L_2}$ is a knot) is the trefoil shown in Figure~\ref{fig:trefoil}.  See Figure~\ref{fig:l5k5} for $\ddddot{L_5}$. 
We will see in Lemma~\ref{lemma:comp}  that all saturated standard hexagonal $r$-mosaics including $\ddddot{L_r}$ contain $r-1$ components.

For $r>3$, we can take the link $\ddddot{L_r}$ and form an alternating knot $\ddddot{A_r}$ from it by removing $r-2$ of the existing (3 crossing) interior tiles and replacing them with  Tile 21, a two crossing tile, in Figure~\ref{fig:tiles}. 
The replacement is equivalent to smoothing one of the three crossings on each of the original tiles.   Any time we smooth a crossing between 2 different components we reduce the number of components of the link by one.  There is one component of $\ddddot{L_r}$ which crosses itself three times in the most central tile (the $0^{th}$ corona), call it $c_1$. We see $c_1$ labeled in $\ddddot{L_5}$ depicted on the left in Figure~\ref{fig:l5k5}.  The horizontal interior arc 
of $c_1$ that passes through that center-most tile crosses all of the other link components on central tiles (two components twice and the rest four times).  We can then choose one of these central crossing tiles for each of the other $r-2$ components and smooth them to change $\ddddot{L_r}$ from an $r-1$ component link into a knot.      Smoothing a crossing in an alternating link results in another alternating link so the end result of this process is an alternating knot.  We see the process of smoothing 3 crossings to turn $\ddddot{L_5}$ into $\ddddot{A_5}$ as described above in Figure~\ref{fig:l5k5} with the smoothed tiles marked $s_1$, $s_2$, and $s_3$.
 Call the alternating knot $\ddddot{A_r}$. $\ddddot{A_r}$, of course, need not be uniquely defined as there are many other ways to choose which crossings to smooth, but any of the choices will work as long as we restrict to the central tiles as this prevents nugatory crossings which might be created by smoothing in the penultimate corona.

\begin{figure}[tpb]  \centering
	\setlength{\unitlength}{0.1\textwidth}
	\begin{picture}(9,3)
	
	\put(0,-.057){\includegraphics[width=.32\textwidth]{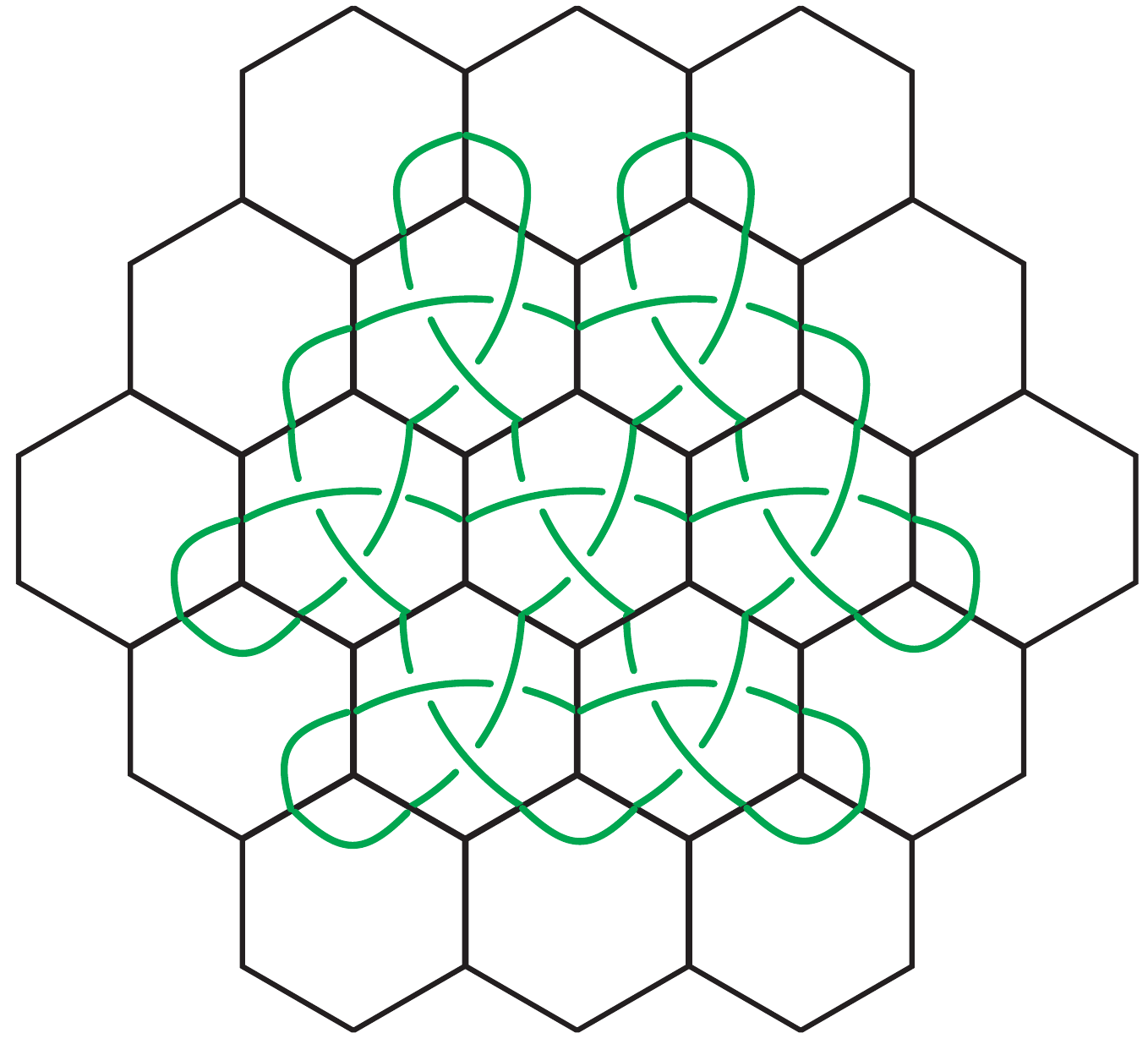}}
	\put(3.3,-.057){\includegraphics[width=.32\textwidth]{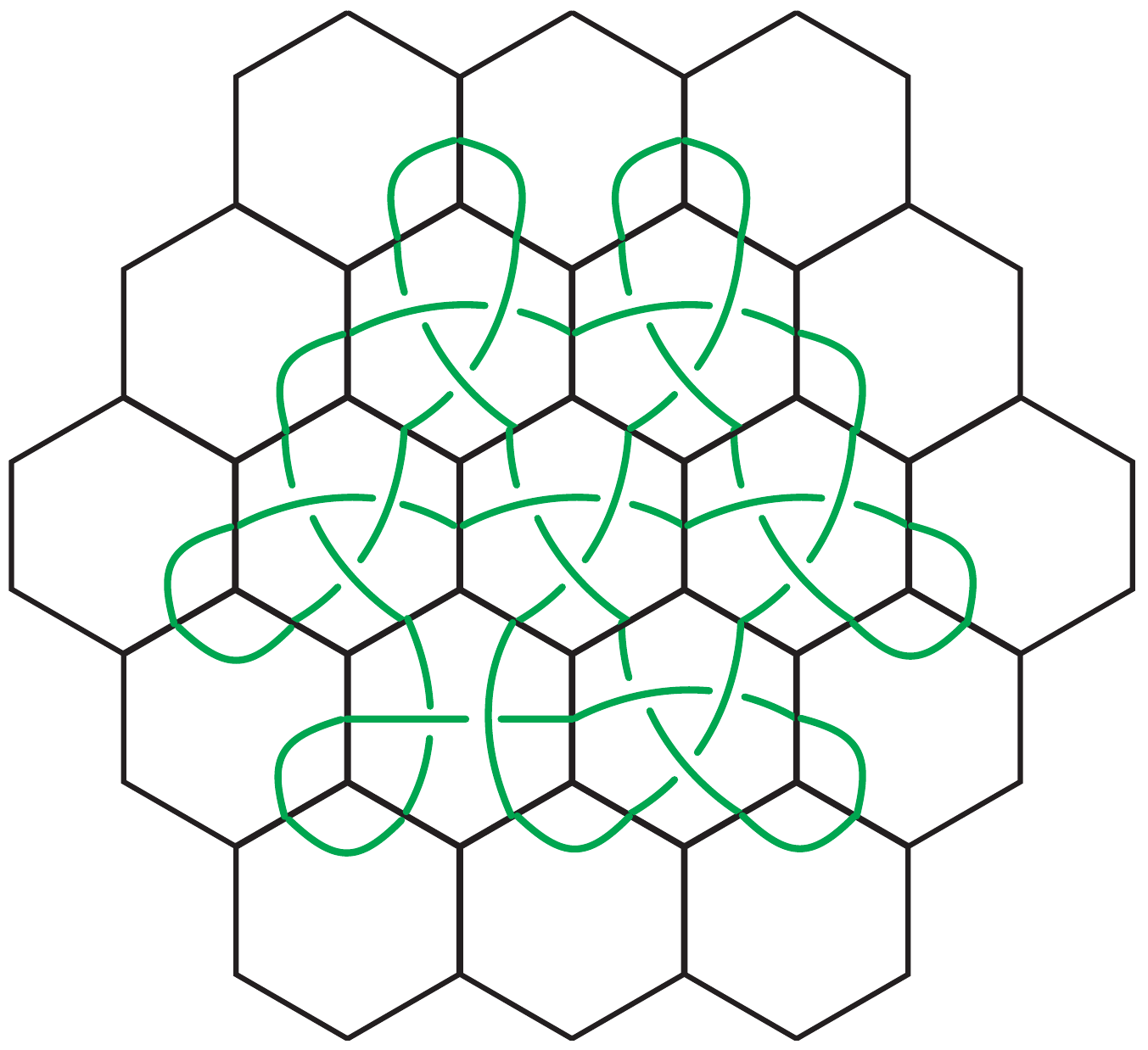}}
	\put(6.6,-.057){\includegraphics[width=.32\textwidth]{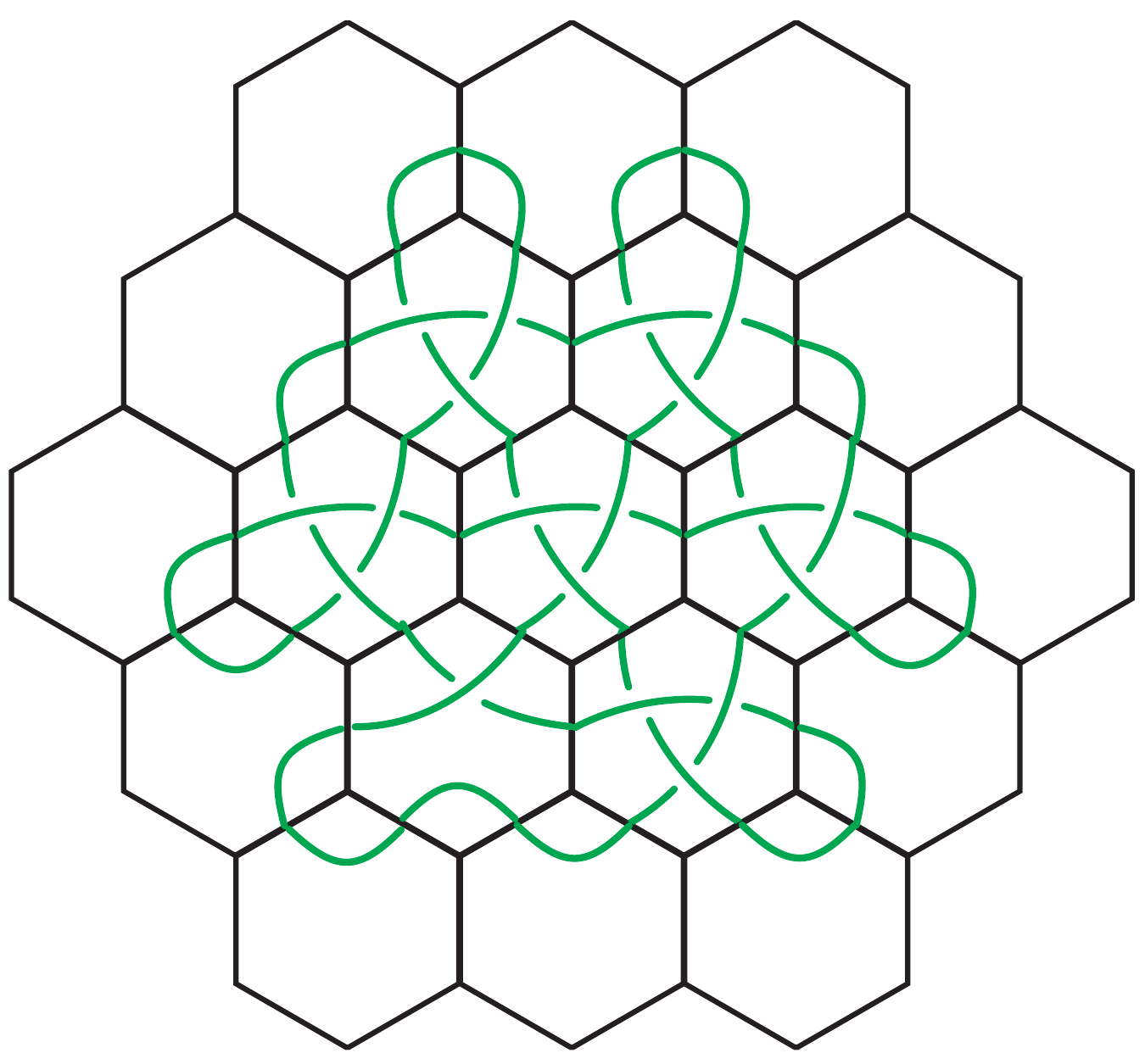}}
	\put(1.25,.56){$c$}
	\put(6.9,2.5){$\ddddot{A_3}$}
	\put(0.3,2.5){$\ddddot{L_3}$}

	\end{picture}
	\caption{Here $\ddddot{L_3}$ is on the left and $\ddddot{A_3}$ on the right.  The central picture is the result of smoothing a crossing of $\ddddot{L_3}$, but before removing a nugatory crossing to form $\ddddot{A_3}$.  Note $\ddddot{L_3}$ has 2 components, but a reduced alternating knot cannot be obtained by smoothing a single crossing so $\ddddot{A_3}$ has 2 fewer crossings that $\ddddot{L_3}$.  Here we smoothed $c$, the bottom left crossing in $\ddddot{L_3}$, to get the middle picture and then did a type I move (or equivalently a second smoothing on the same tile) to get $\ddddot{A_3}$ on the right.
	Because of all the symmetry, it is easy to check that $c$ is not unique and that no crossing in $\ddddot{L_3}$ can be smoothed to immediately result in a reduced alternating knot.}
	\label{fig:l3k3}
\end{figure}
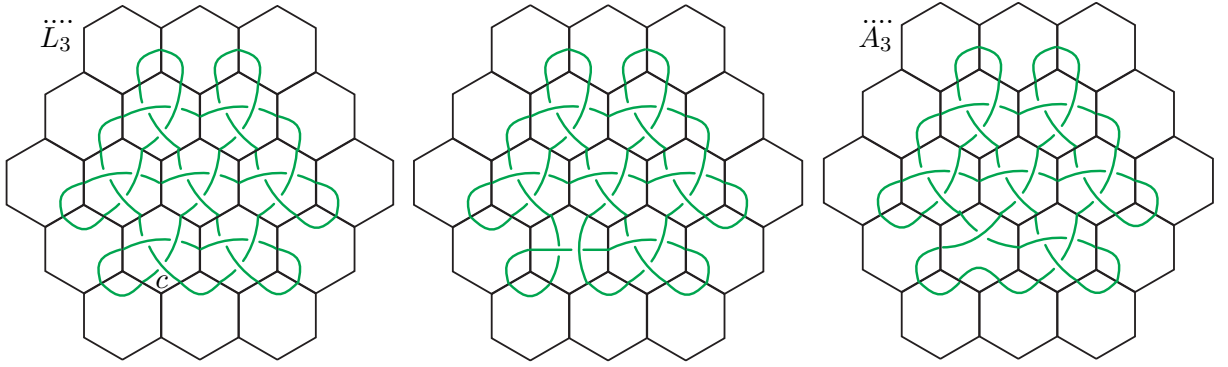

Note that this process works for $r>3$ and works trivially for $r=2$, but does not work for $r=3$.  The problem for $r=3$ is that there are only 2 components in $\ddddot{L_3}$ and all of the arcs on the central tiles belong to a single component (since the entire set of central tiles consists of only a single tile) as you can see in Figure~\ref{fig:l3k3}. All of the crossings between distinct components occur in the penultimate corona and up to symmetry they break down into two equivalence classes. Because smoothing must be done between distinct components and thus in the penultimate corona, no matter which of the two types of crossings we pick the tile substitution trick above results in a knot with a nugatory crossing instead of a reduced alternating knot as it did in all the other cases where the substitution could be done on one of the central tiles. As a result  $\ddddot{A_3}$ has crossing number 19, one less than would have been predicted simply by taking a saturated mosaic and smoothing $r-2$ crossings to obtain a knot.   This problem is special to  $\ddddot{A_3}$ and does not cause problems for $\hat{A_3}$ or $\widehat{A_3}$ as we'll see below.

\begin{figure}[tpb]  \centering
	\setlength{\unitlength}{0.1\textwidth}
	\begin{picture}(9,3)
	
	\put(0,-.057){\includegraphics[width=.32\textwidth]{a3dotsai.pdf}}
	\put(3.3,-.057){\includegraphics[width=.32\textwidth]{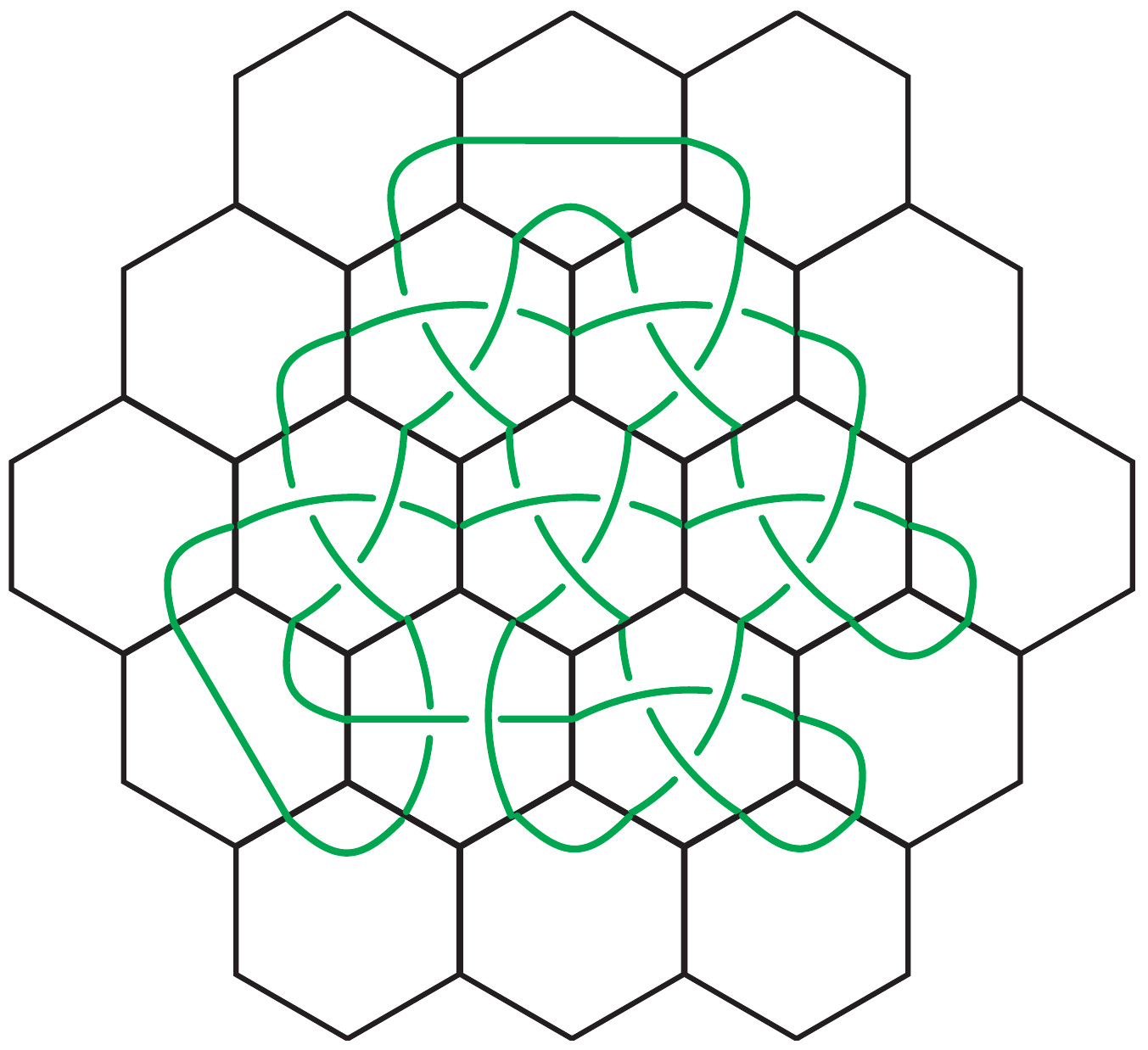}}
	\put(6.6,-.057){\includegraphics[width=.32\textwidth]{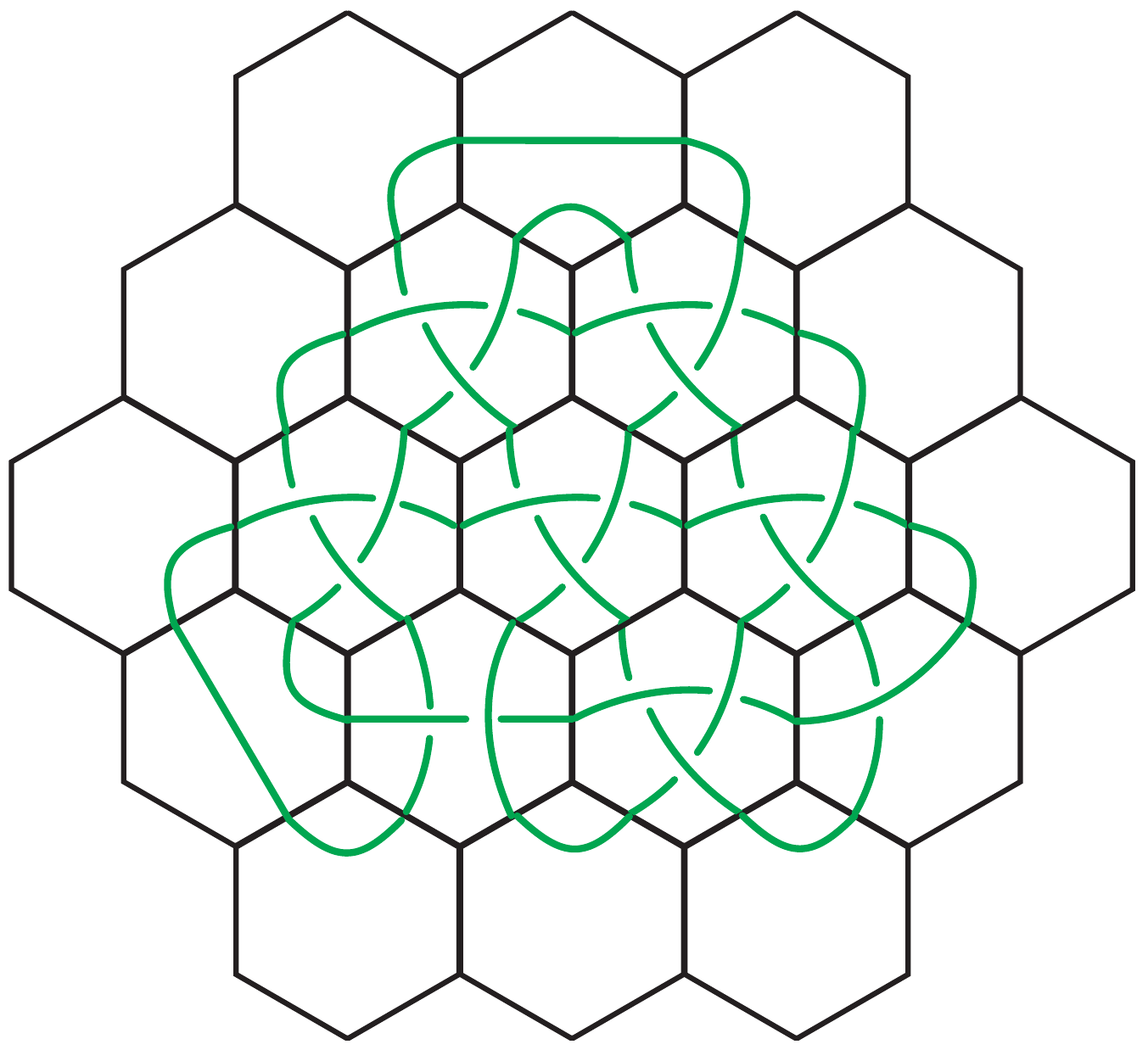}}
	\put(1.25,.56){$c$}
	\put(6.9,2.5){$\widehat{A_3}$}
	\put(3.6,2.5){$\hat{A_3}$}
	\put(0.3,2.5){$\ddddot{A_3}$}

	\end{picture}
	\caption{We see a knot with maximal crossing number on a 3-mosaic in each of the settings. Standard  $\ddddot{A_3}$ is on the left,  semi-enhanced $\hat{A_3}$ in the middle and enhanced $\widehat{A_3}$ on the right.}
	\label{fig:three}
\end{figure}

In the semi-enhanced setting, we can form $\hat{L_r}$ from $\ddddot{L_r}$ by swapping tiles of Type 6 for tiles of Type 5 in the boundary corona whenever this lowers the number of components and then we change any necessary crossings to make the new link alternating.  We   then form $\hat{A_r}$ from $\hat{L_r}$ (for $r>3$) by smoothing crossings on central tiles to get a reduced alternating knot.  We see the process of forming $\hat{L_7}$ from $\ddddot{L_7}$ in Figure~\ref{fig:tile7}.  For $r=3$ we see $\hat{A_3}$ depicted in the middle in Figure~\ref{fig:three}.

\medskip
\emph{The definitions of $\widehat{L_r}$ and $\widehat{A_r}$ on enhanced hexagonal $r$-mosaics:}
\medskip

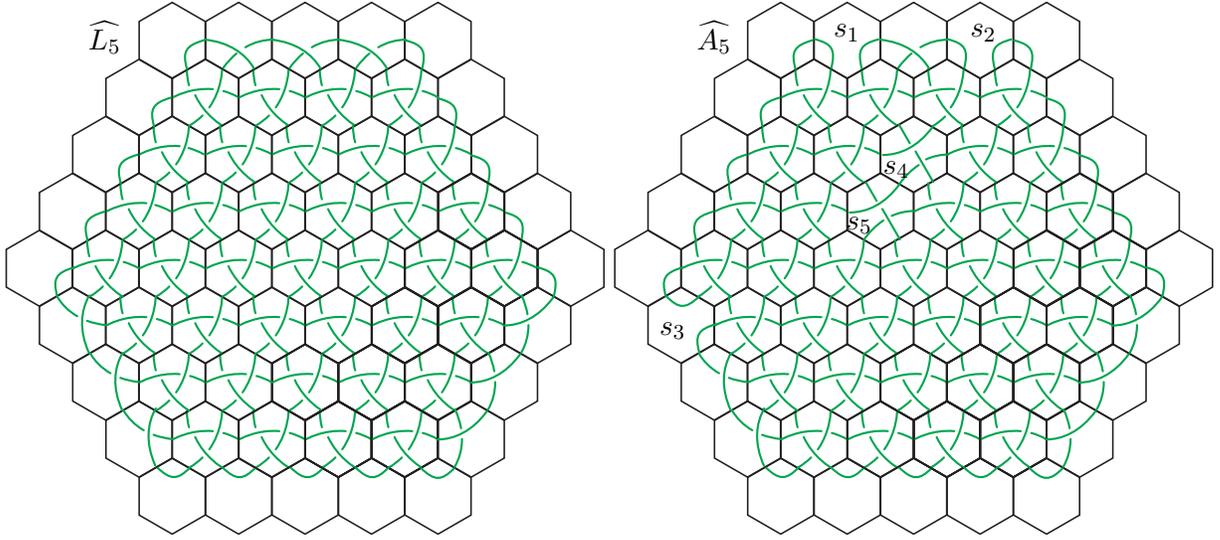
\begin{figure}[tpb]  \centering
	\setlength{\unitlength}{0.1\textwidth}
	\begin{picture}(9,4.5)
	\put(0.1,-.057){\includegraphics[width=.49\textwidth]{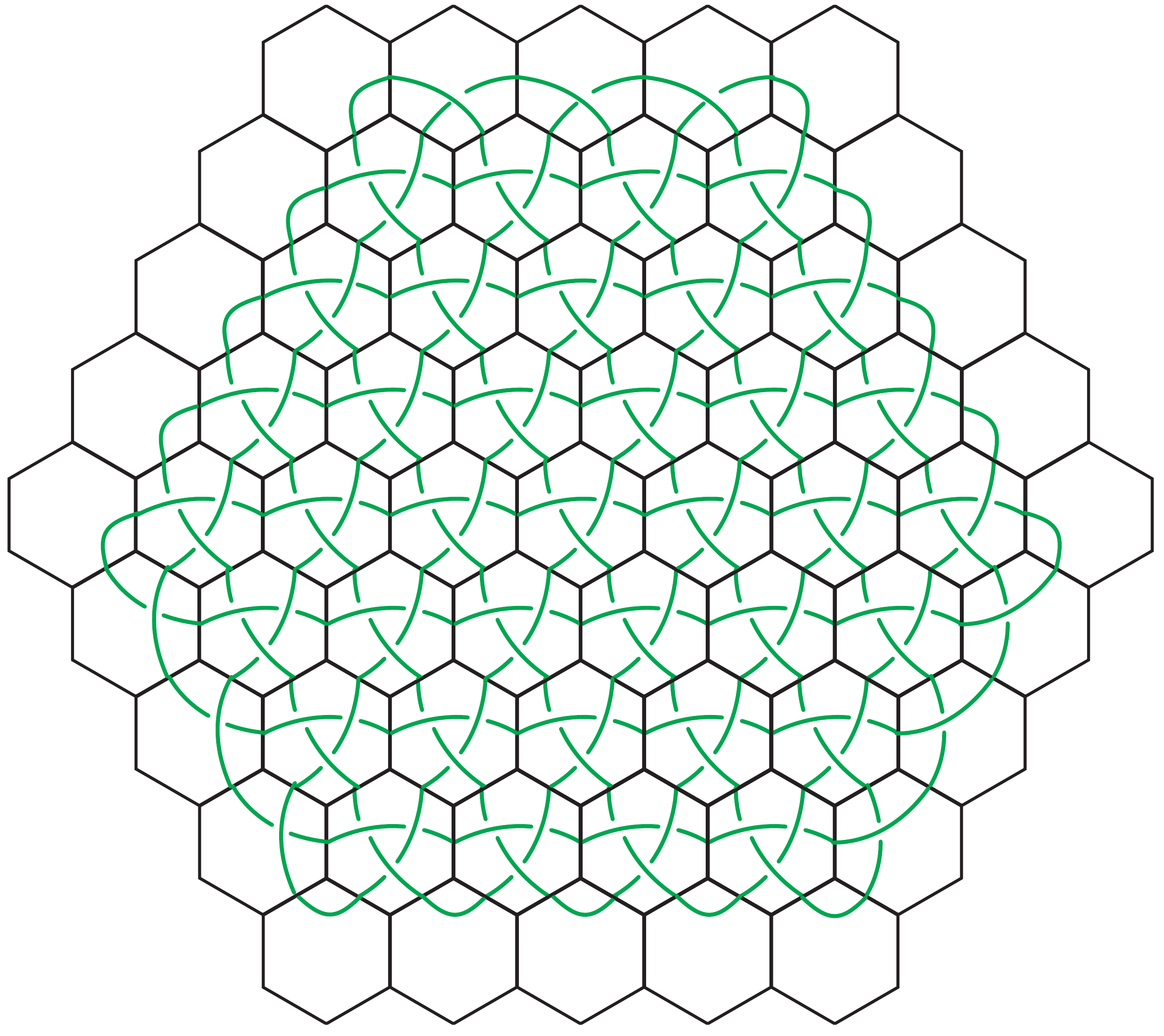}}
	\put(5,-.057){\includegraphics[width=.49\textwidth]{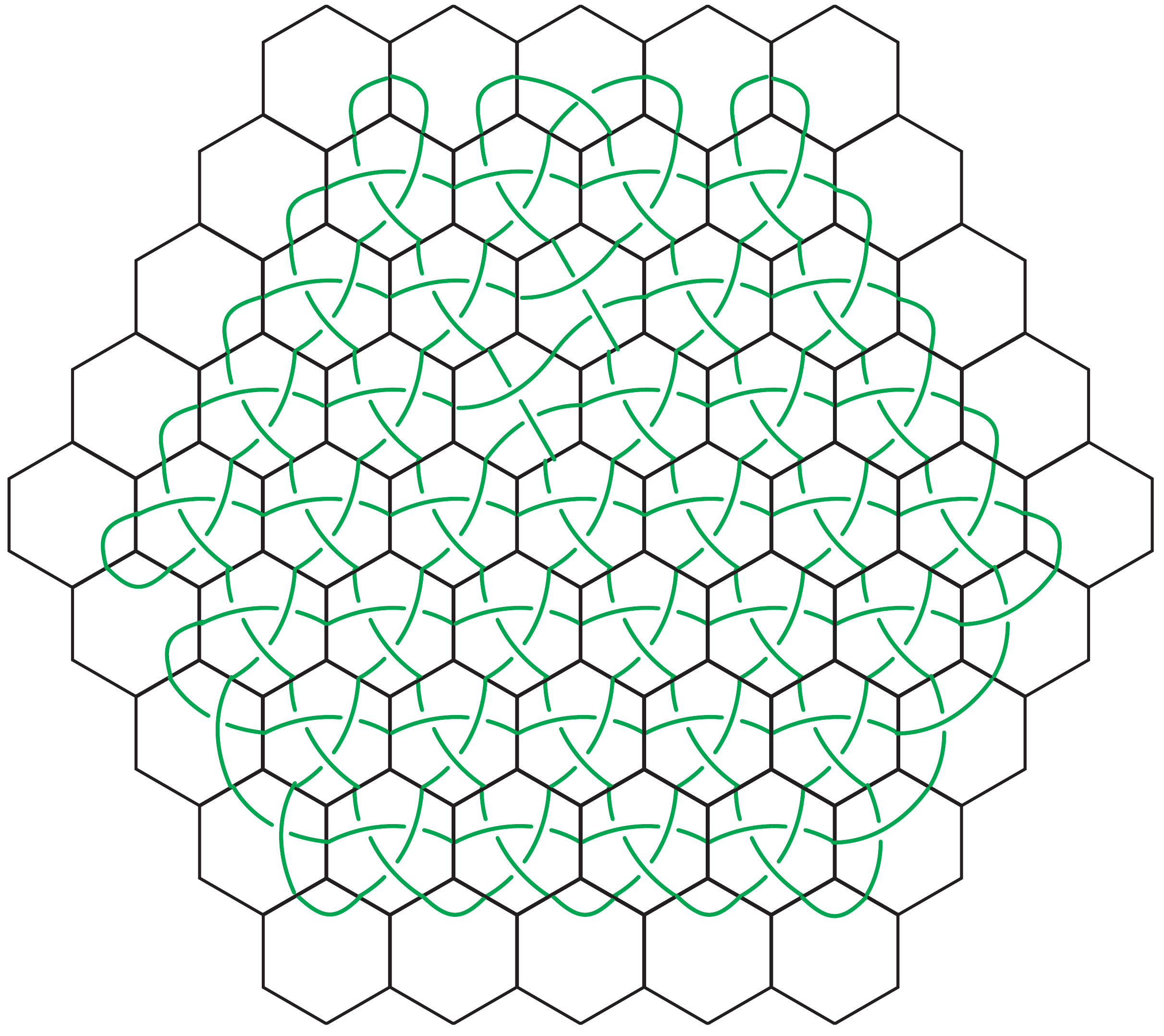}}
		\put(0.8,3.9){$\widehat{L_5}$}
	\put(5.7,3.9){$\widehat{A_5}$}

	\put(6.8,4){$s_1$}
	\put(7.9,4){$s_2$}
	\put(5.4,1.6){$s_3$}
	\put(7.2,2.9){$s_4$}
	\put(6.9,2.45){$s_5$}

	\end{picture}
	
\label{fig:l5a5Sat}
	\caption{Here on the left $\widehat{L_5}$ in the enhanced hexagonal setting is converted into an alternating knot $\widehat{A_5}$ on the right
by smoothing crossings. Smoothings $s_1$, $s_2$, and $s_3$ combine bigons with other components by replacing a boundary crossing tile by a boundary non-crossing tile and smoothings $s_4$ and $s_5$ replace three crossing tiles with two crossing tiles giving the final necessary smoothings on central tiles.  The same general strategy is used for $\widehat{A_r}$ for all $r \geq 4$.}
\label{fig:l5satk5sat}
\end{figure}

Now we define $\widehat{L_r}$ and $\widehat{A_r}$ in the enhanced hexagonal setting.  The interior tiles of  $\widehat{L_r}$ are saturated and identical to the interior tiles for $\ddddot{L_r}$.  The boundary tiles then are filled in  alternating sides of the boundary corona filled with crossing tiles and sides without crossings like in $\widehat{L_5}$ the board on the left in Figure~\ref{fig:l5satk5sat}.  $\widehat{A_r}$ is then formed from $\widehat{L_r}$ by smoothing crossings as before, depicted in Figure~\ref{fig:l5satk5sat}.  As we will see below, $\widehat{L_r}$ contains $r+1$ components  for $r>2$ and $r$ smoothings will combine these components into a single knot.
As before we can form the alternating knot $\widehat{A_r}$  by replacing interior tiles, but now since there are crossings in boundary tiles we can also choose to smooth crossings of distinct components of $\widehat{L_r}$ that occur in the boundary tiles if we prefer.

$\widehat{L_r}$ is a knot when $r=2$  so $\widehat{L_2} = \widehat{A_2}$, so we need only describe how to form  to $\widehat{A_r}$ for $r > 2$.
Rotate $\widehat{L_r}$ if necessary so that the top edge of boundary tiles contains crossings as it does in the figures in this paper. 
To form $\widehat{A_r}$ we note that $\widehat{L_r}$ contains three components running between diametrically opposed corners of the boundary corona.  We call these three components bigons since they have only two interior arcs.  $\widehat{L_r}$ also always contains one component with three interior arcs which we call a triangle. If $r$ is odd the triangular component  hits the center tile of the three sides of the boundary corona which have no crossings as seen in $\widehat{L_5}$ depicted in Figure~\ref{fig:l5satk5sat}.  If $r$ is even  it hits the two most central crossing tiles on each of three boundary edges that do contain crossings as seen in $\widehat{L_4}$ in Figure~\ref{fig:l4}.  For both even and odd boards  the remaining components consist of $r-3$ hexagons in the sense that they have 6 interior arcs, two horizontal ones, two with positive slope and two with negative slope.

Call the bigon component which hits $T_{1,1}$, the top left tile, $c_1$.  As we move left to right across the top edge of tiles call the second component we come to $c_2$ the third component $c_3$ and so on until we get to tile $T_{1,r}$ in the top right corner which intersects the second bigon which we will call $c_j$ (the argument does not depend on the exact value of $j$, but to be precise if $r=2k-1$ then $j=k$ if $r=2k$ then $j=k+1$).  Call the third bigon $c_{j+1}$.  $c_{j+1}$ is horizontal, does not intersect the top edge and runs from $T_{r,1}$ on the far left to tile $T_{r,2r-1}$ on the far right.  Note that in the top edge $c_1$ shares a crossing with $c_2$, $c_2$ shares a crossing with $c_3$ and so on until $c_{j-1}$ shares a crossing with $c_j$.   Smooth the boundary crossing between $c_1$ and $c_2$ to combine these into a single link component which we will call $c_{12}$ (see smoothing $s_1$ in Figure~\ref{fig:l5satk5sat}).  Now smooth the crossing that had been between $c_2$ and $c_3$ to form component $c_{123}$.  Continue until we have smoothed $j-1$ crossings and formed $c_{123 \dots j-1 j}$ (smoothing $s_2$ in Figure~\ref{fig:l5satk5sat}).

The third bigon $c_{j+1}$  crossed $c_2$ in the boundary tile $T_{r+1,1}$ so now it crosses $c_{123 \dots j-1 j}$ in that tile.  Smooth that crossing to form $c_{123 \dots j j+1}$ again decreasing the number of components by one (smoothing $s_3$ in Figure~\ref{fig:l5satk5sat}).  
If $r=3$ all we have left besides $c_{123}$ is a triangular component.  In that case, we smooth a crossing between these two components to form $\widehat{A_3}$ as in Figure~\ref{fig:three}. So now we turn our attention to $r>3$.
If $r$ is even the triangle was one of the smoothed components  combined with $c_1$ and
 all the remaining unsmoothed components are  $\frac{r-2}{2}$ 
  hexagons that were disjoint from the boundary tiles of $\widehat{L_r}$ containing crossings.  
 If $r$ is odd the remaining components consist of one triangle and  
 $\frac{r-3}{2}$ 
 hexagons  that were disjoint from the boundary tiles of $\widehat{L_r}$ containing crossings. 
 On the odd board ($r>3$), the triangle has a horizontal interior arc in the top half of the mosaic that crossed $c_1$ twice on central tiles.
 In both the even and odd cases, each of the remaining hexagonal components has two horizontal interior arcs in the top half of the board.  The lower of these two arcs crossed $c_1$ on a central tile (twice, of course).  This means that they all now cross $c_{123 \dots j j+1}$ on a central tile.  Smooth one of those crossings for each component to form $\widehat{A_r}$ (smoothings $s_4$ and $s_5$ in Figure~\ref{fig:l5satk5sat}).  Again note that this process forms a reduced, alternating knot smoothing as few crossings as possible, but our choices were far from unique. Many other crossings could have been smoothed to create a reduced, alternating knot with the same crossing number.  $\widehat{L_r}$ had $r+1$ components so we note that the crossing number of $\widehat{A_r}$ is exactly $r$ lower than the crossing number of $\widehat{L_r}$.  The crossing number of $\widehat{L_r}$ (enhanced) is the crossing number of $\ddddot{L_r}$ (standard) plus $3(r-2)$ because they each have the same number of crossings on the interior tiles, but $\widehat{L_r}$ has $3(r-2)$ crossings on its boundary tiles.

\section{Useful Lemmas}
In this section we establish some useful results for hexagonal $r$-mosaics.  The equivalent lemmas are either clear or well known in the rectangular setting.

\label{sec:lud}

\begin{lemma}  Let $\ddddot{M_1}$ and $\ddddot{M_2}$ be saturated standard hexagonal $r$-mosaics. 
A rotation of $\ddddot{M_1}$ by 0 or $\frac{\pi}{3}$ together with a combination of crossing changes and planar isotopies of individual tiles is sufficient to get from $\ddddot{M_1}$ to $\ddddot{M_2}$.
\label{lemma:iso}
\end{lemma}

\begin{proof}
All of the connection points of the boundary of the interior tiles are used  by $\ddddot{M_1}$ and $\ddddot{M_2}$ since they are both saturated.  
There are two ways to connect a standard saturated mosaic up through the boundary tiles.  If the boundary tiles of $\ddddot{M_1}$ and $\ddddot{M_2}$ do not already match a rotation of $\ddddot{M_1}$  by $\frac{\pi}{3}$ will cause the two sets of boundary tiles to match.  Now we just need to work on the interior tiles.  All interior tiles of $\ddddot{M_1}$ and $\ddddot{M_2}$ are three crossing tiles (Tiles 24-27 in Figure~\ref{fig:tiles}).    Any two such tiles are homotopic (rel $\partial$) and thus we can get from any one of these types of tiles to any other with crossing changes and planar isotopies of the individual tiles.  
\end{proof}

\begin{lemma}
All saturated standard hexagonal $r$-mosaics have the same number of components.
\end{lemma}

\begin{proof}
Since crossing changes, planar isotopies and rotations do not alter the number of components this follows directly from the previous lemma.
\end{proof}

\begin{lemma}
On a standard hexagonal $r$-mosaic saturated links  (including $\ddddot{L_r}$)
 have $r-1$ components for $r \geq 2$.  The enhanced hexagonal $r$-mosaic $\widehat{L_r}$ contains $r+1$ components for $r>2$ and 1 component if $r=2$. The semi-enhanced hexagonal $r$-mosaic $\hat{L_r}$ contains $\lceil \frac{r}{2} \rceil$ components  for $r>2$ and 1 component if $r=2$.
Each $L_r$ has the smallest number of components of any saturated $r$-mosaic in its category.
 \label{lemma:comp}
\end{lemma}

\begin{figure}[tpb]  \centering
	\setlength{\unitlength}{0.1\textwidth}
	\begin{picture}(9,5)
	\put(-0.2,-.057){\includegraphics[scale=.75]{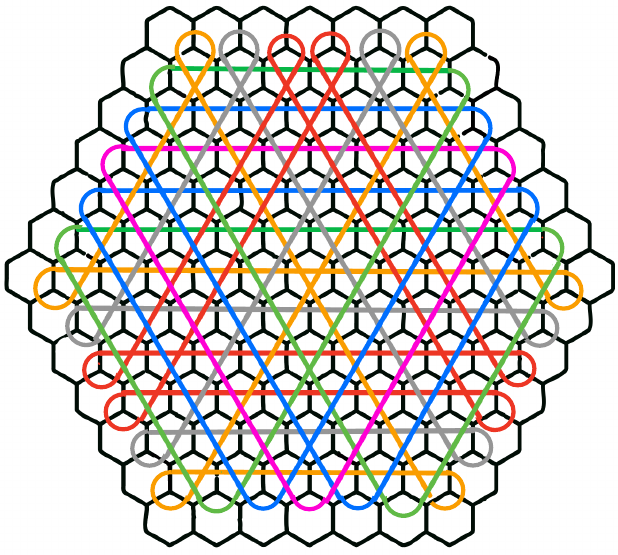}}
	\put(4.5,-.057){\includegraphics[scale=.75]{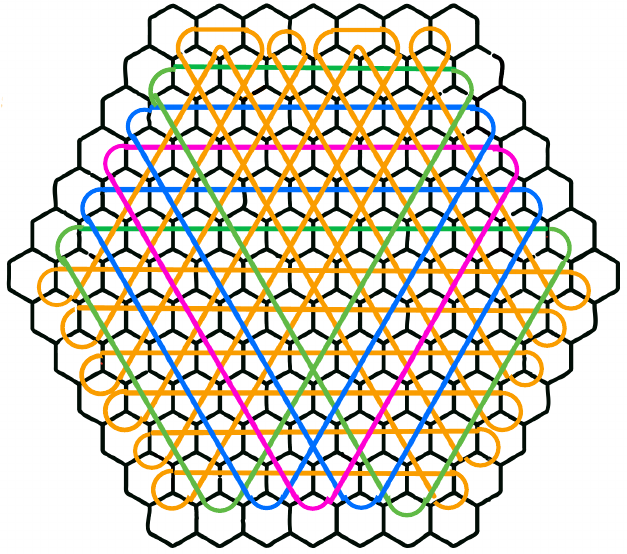}}
	\put(5.25,3.98){$\hat{L_7}$}
	\put(.55,3.98){$\ddddot{L_7}$}
	\put(.99,3.98){$c_1$}
	\put(1.43,3.98){$c_2$}
	\put(2.5,3.98){$c_2$}
	\put(2.08,3.98){$c_3$}
	\end{picture}
	\caption{Here the shadow of $\ddddot{L_7}$ is on the left and $\hat{L_7}$ on the right. $\ddddot{L_7}$ has 6 components as expected by \cite{jmm}.  $\hat{L_7}$ shows that by allowing Tile 6 from Figure~\ref{fig:tiles} in the boundary tiles we can get a saturated semi-enhanced hexagonal 7-mosaic with only 4 components combining $c_1$ and $c_2$ in one boundary tile and $c_2$ and $c_3$ in another.  We do not worry about over and under crossings here since the figure is crowded and changing crossings does not have any impact on the number of components in a link.} 
	\label{fig:tile7}
	\label{fig:LPrime}

\end{figure}
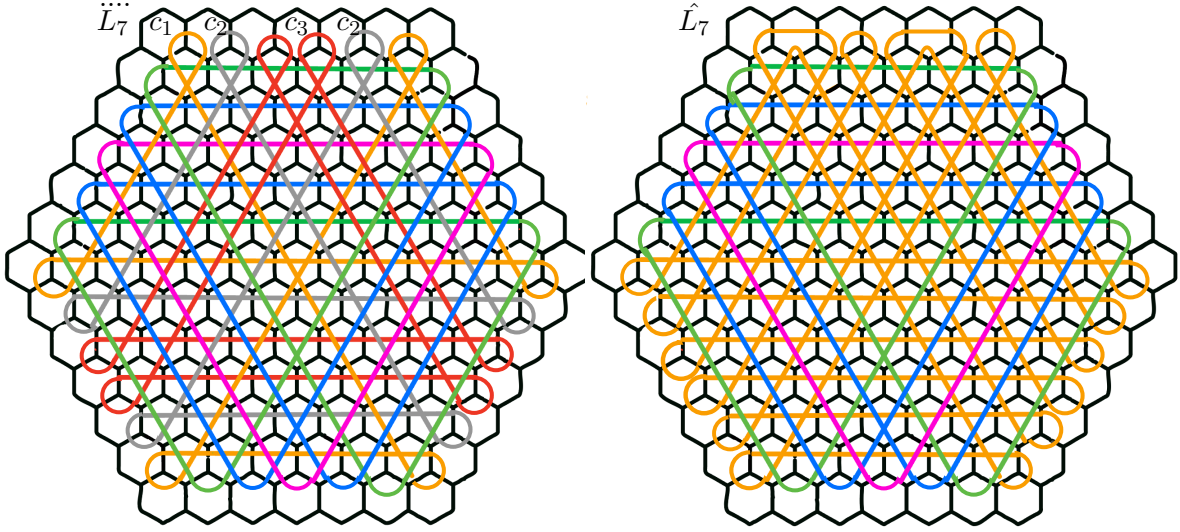

\begin{proof}
The standard case is a special case of the Theorem 3.1 in \cite{jmm}.   In an enhanced 2-mosaic there cannot be any crossings in the boundary tiles because all boundary tiles in that case are corner tiles so $\widehat{L_2}=\ddddot{L_2}=\hat{L_2}=\widehat{K_2}=\ddddot{K_2}=\hat{K_2}$.  
As mentioned previously for $r > 2$ $\widehat{L_r}$ contains 3 bigon components, one triangle, and $r-3$ hexagonal components.  

Note that $\hat{L_r}$  shows Theorem 3.1 in \cite{jmm} would not hold for semi-enhanced hexagonal mosaics.  For example, $\hat{L_7}$ has only 4 components as shown by Figure~\ref{fig:LPrime}.  The difference between standard hexagonal mosaics and semi-enhanced is that
 for standard hexagonal mosaics  such as $\ddddot{L_r}$ each of the boundary tiles, which used four connection points had to be a Type 5 tile, but in the semi-enhanced setting any of these tiles could be replaced by Type 6 tiles.  
Swapping the tiles is equivalent to doing a band sum between two arcs of a  link (or knot).  The band sum will decrease the number of components if and only if it is between distinct components of the link.   Under those conditions,  
 the bands allow us to turn $\ddddot{L_r}$ into a saturated link   $\hat{L_r}$ with 
 only
  $\lceil \frac{r}{2} \rceil$ components.  Every saturated semi-enhanced hexagonal mosaic can be formed from a saturated standard hexagonal mosaic swapping out some Type 5 boundary tiles for Type 6 boundary tiles and thus one can show that $\lceil \frac{r}{2} \rceil$ is indeed the smallest number of components on any  saturated semi-enhanced hexagonal $r$-mosaic.

  We see the shadows of $\ddddot{L_7}$ and $\hat{L_7}$ in Figure~\ref{fig:tile7}. In $\ddddot{L_7}$  we have the expected 6 components.  Components $c_1$ and $c_2$ enter tile $T_{1,2}$ and $c_2$ and $c_3$ enter tile $T_{1,5}$ allowing us to insert Type 6 boundary tiles reducing the number of components in the semi-enhanced case from 6 to 4.   
    
  The general formula holds because for $r=4$ or 5 the smoothings allow us to combine 2 components dropping the total number in the link by one, for $r=6$ or 7 we can combine up to 3 components dropping the number by two, and each time $r$ increases by two we get one more component that can be combined using Type 6 boundary tiles.  
  $\ddddot{L}$ has exactly $r-1$ components  and $\lfloor \frac{r}{2} \rfloor$ components which hit each of the edges of the boundary corona containing the type 5 tiles.  Thus we can drop by as much as $\lfloor \frac{r}{2} \rfloor -1$ components giving a lower bound of $r-1- (\lfloor \frac{r}{2} \rfloor -1)) = r - \lfloor \frac{r}{2} \rfloor =  \lceil \frac{r}{2} \rceil$.  We, however, achieve this bound exactly for $\hat{L_r}$ and thus it achieves the smallest number of components in a saturated semi-enhanced $r$-mosaic.

\end{proof}

\begin{claim}
Let $\ddddot{s_r}$ be the number of crossings of a  saturated standard hexagonal $r$-mosaic and $\hat{s_r}$ be number of crossings of a  saturated semi-enhanced hexagonal $r$-mosaic. We observe that $\ddddot{s_r} = \hat{s_r}  =  9r^2-27r+21$ for $r \geq 2$.    
\end{claim}

\begin{proof}

In a saturated standard  (or semi-enhanced) hexagonal mosaic all crossings are on interior tiles and there are 3 crossings on each interior tile, so $\ddddot{s_r}= \hat{s_r}$ is 3 times the well known \emph{centered hexagonal number} which tells how many hexagons are in a centered hexagonal grid or radius $r-1$ so $\ddddot{s_r} = (3*(3r^2-9r+7)) = 9r^2-27r+21$.  

\end{proof}

\begin{claim}
Hexagonal $r$-mosaics $\ddddot{L_r}$ (standard) and $\hat{L_r}$ (semi-enhanced) have crossing number $\ddddot{s_r}  = \hat{s_r}  = 9r^2-27r+21$ for $r \geq 2$.
\label{claimlr}
\end{claim}

This is clear since there $\ddddot{L_r}$ and $\hat{L_r}$ have 3 crossings on each interior tile, none on the boundary tiles, contain no nugatory crossings, and are alternating.

	\begin{lemma}  An enhanced hexagonal $r$-mosaic cannot have crossings on adjacent sides of the (hexagonal) boundary corona without an endpoint of an arc of the complement falling between those crossing tiles.  
	\label{lemma:boundarycrossings}
	\end{lemma}
	 
	 \begin{proof}
	 	Examine a corner tile of the boundary corona.  Without loss of generality we let it be tile $T_{1,1}$.  Every corner tile has 3 potential connection points since it has 3 sides that are disjoint from other tiles of the mosaic and three that are not.  Thus at most 2 of these three connection points can be used by a link
	and one is unused.  If the unused connection point intersects $T_{2,2}$, the adjacent interior tile,  then we have an arc of the complement starting there as the lemma claims.  If not then one of the connection points of one of the adjacent boundary tiles must be missed.  Without loss of generality let it miss the connection point between $T_{1,1}$ and $T_{2,1}$.   This then means that $T_{2,1}$ cannot use all 4 of its possible connection points thus it contains only one arc and no crossings.  If that arc misses one of the two interior connection points of $T_{2,1}$ we again have an arc of the complement there satisfying the lemma.  If not then it must miss the connection point between $T_{2,1}$ and $T_{3,1}$.  The same argument is now made for $T_{3,1}$ as we just made for $T_{2,1}$.  This pattern repeats until we either miss a connection point with the interior yielding an arc of the complement or we reach the next corner tile ($T_{r,1})$ without any crossings in the boundary tiles. 
\end{proof}

\begin{claim}
The enhanced hexagonal $r$-mosaic $\widehat{L_r}$ has crossing number $  9r^2-24r+15$ for $r \geq 2$ and no enhanced hexagonal $r$-mosaic can have more crossings than this.
\label{claimlrenh}
\end{claim}

\begin{proof} By Lemma~\ref{lemma:boundarycrossings}  
a saturated enhanced hexagonal mosaic can have crossings on up to three of the edges of boundary tiles and no more.    Since each of these have $r-2$ possible crossing tiles we get	$\ddddot{s_r} +3(r-2)= 9r^2-24r+15$, which matches $\widehat{L_r}$.
\end{proof}

\begin{claim}
On a standard hexagonal $r$-mosaic, the knot $\ddddot{A_r}$ has crossing number $\ddddot{s_r} - (r-2) = 9r^2-28r+23$ for $r=2, r > 3$. It has crossing number  $19= 9r^2-28r+22$ when $r=3$ because there is no way to smooth the crossings of $\ddddot{L_3}$ without creating one nugatory crossing.  For the semi-enhanced setting, $\hat{A_r}$ has crossing number $\ddddot{s_r} - (\lceil \frac{r}{2} \rceil-1)
=  9r^2-27r+22 - \lceil \frac{r}{2} \rceil$ for $r > 2$.  Of course, $\hat{A_2}=\ddddot{A_2}$.  \label{claimar}
\end{claim}

\begin{proof}
	 
	 This is clear since for $r \neq 3$ $\ddddot{A_r}$ is formed from $\ddddot{L_r}$ by smoothing $r-2$ crossings and $\hat{A_r}$ is formed from $\hat{L_r}$ by smoothing $\lceil \frac{r}{2} \rceil-1$ crossings.  For $r=2$ and 3, only one component of  $\ddddot{L_r}$ enters the Type 5 boundary tiles so swapping Type 5 and 6 tiles in the boundary does not decrease the number of components.

\end{proof}

\begin{claim}
On an enhanced hexagonal $r$-mosaic the knot $\widehat{A_r}$ has crossing number $\ddddot{s_r} +3(r-2)-r= 9r^2-25r+15$ for $ r > 2$ and crossing number $3$ for $r=2$.
\label{claimarenh}
\end{claim}

This is clear since $\widehat{A_2} = \ddddot{A_2}$ is a trefoil and for $r > 2$  $\widehat{A_r}$ is formed from $\widehat{L_r}$ by smoothing $r$ crossings  yielding a reduced alternating knot.

\section{Bounding Crossing Number}

\label{sec:cn}

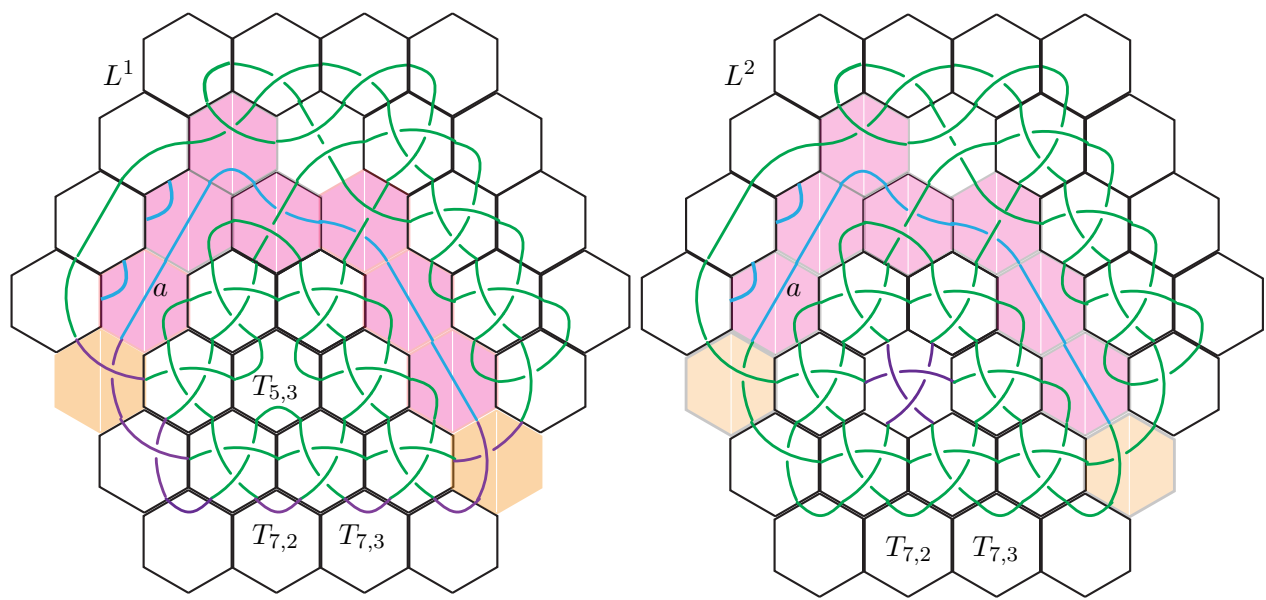
\begin{figure}[tpb]  \centering
	\setlength{\unitlength}{0.1\textwidth}
	\begin{picture}(10,5)
	\put(-0.1,-.057){\includegraphics[width=.51\textwidth]{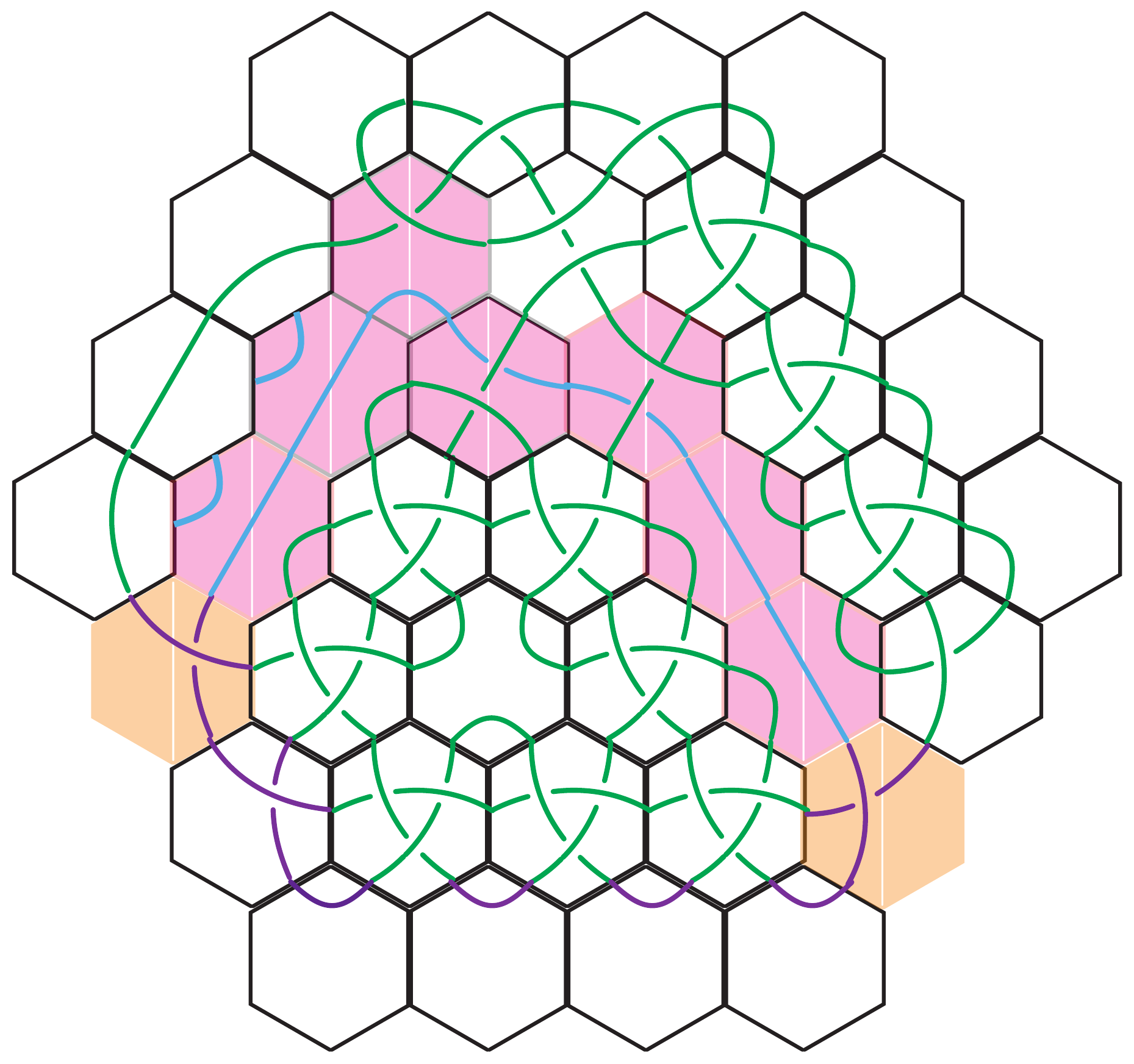}}
	\put(5,-0.06)  {\includegraphics[width=.51\textwidth]{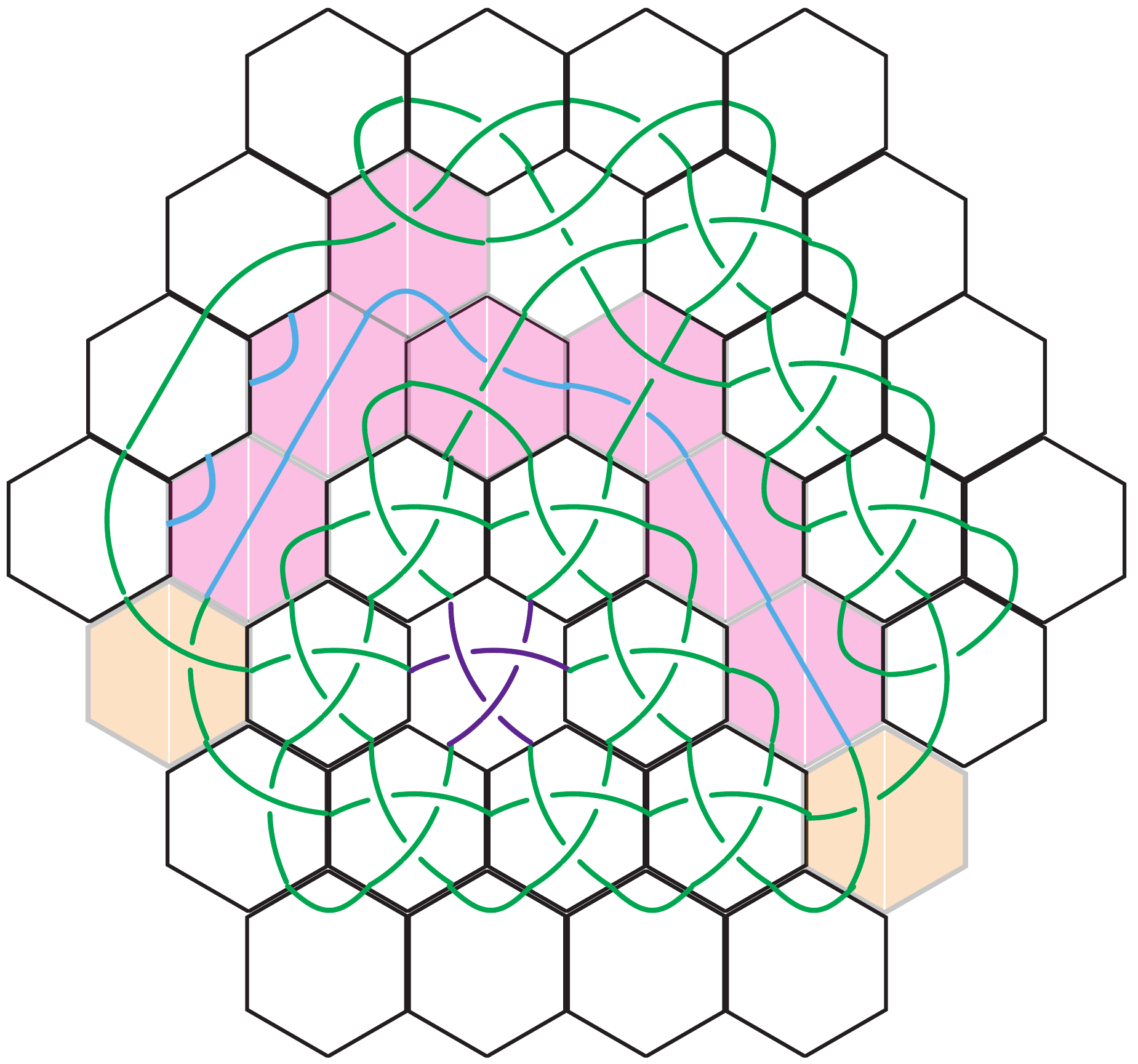}}
	\put(.7,4.1){$L^1$}
	\put(5.7,4.1){$L^2$}
  
	\put(1.9,1.6){$T_{5,3}$} 
	
	\put(1.1,2.4){$a$}
	\put(6.2,2.4){$a$}
	
\put(1.9,0.4){$T_{7,2}$}
	\put(7.0,0.3){$T_{7,2}$}
\put(2.6,0.4){$T_{7,3}$}
	\put(7.7,0.3){$T_{7,3}$}
	
	\end{picture}
	\caption{On the left $L^1$ adds $a$ to the interior arcs of $K^1$ from Figure~\ref{fig:KA}  leaving the boundary tiles in $\tilde{I}$ unchanged and maximizing the number of crossings in the boundary tiles $\tilde{O}$.  The green arcs are the same in the two figures.  The purple is changed.  On the right we form $L^2$ from $L^1$ by replacing and tiles outside $a$ in $O$ which are not already saturated with saturated tiles.  Again the changes are depicted in purple.}
	\label{fig:LA1}
	\label{fig:LA2}
\end{figure}

\begin{figure}[tpb]  \centering
	\setlength{\unitlength}{0.1\textwidth}
	\begin{picture}(10,5)
	\put(-0.1,0){\includegraphics[width=.51\textwidth]{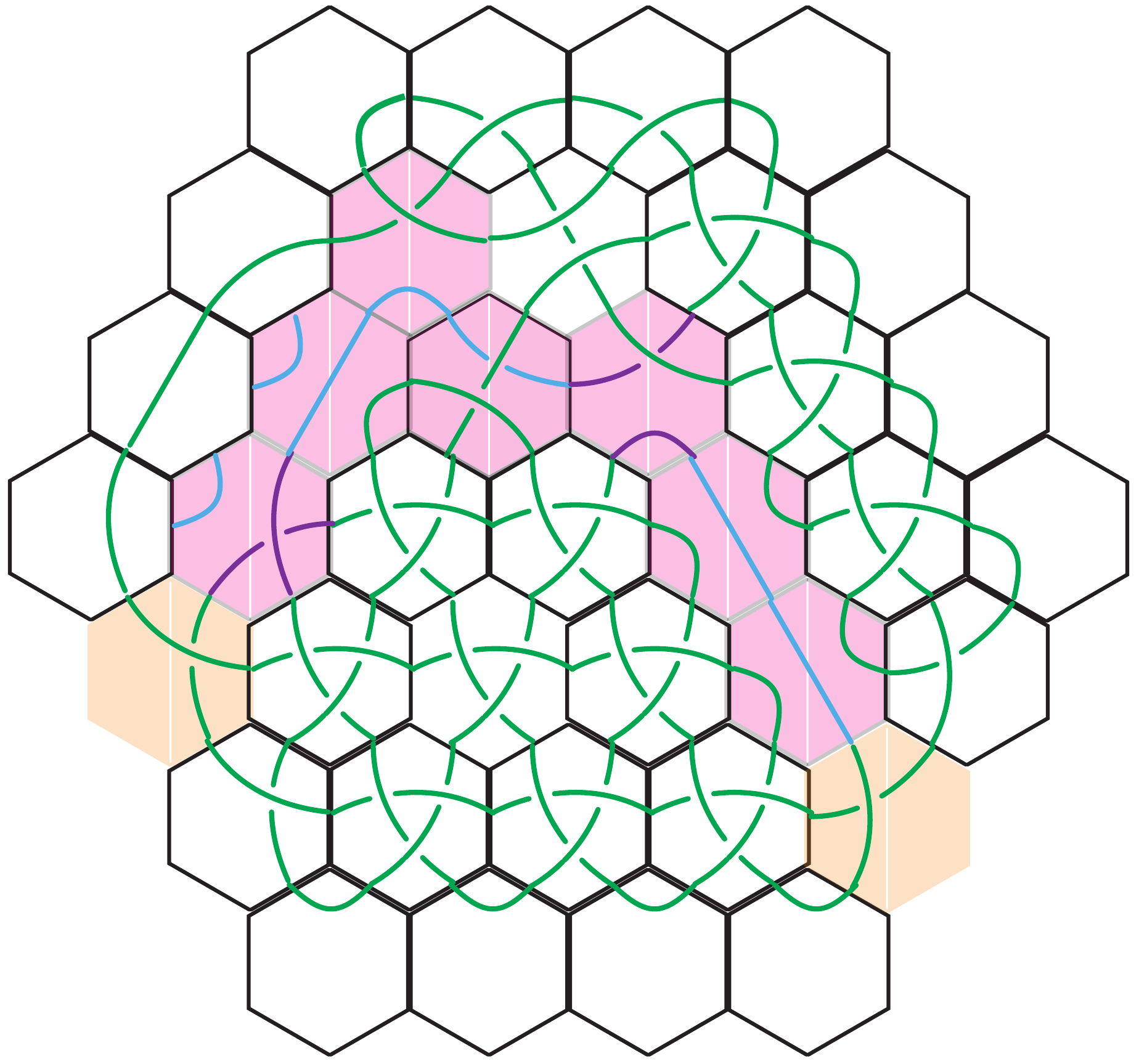}}
	\put(5.2,-.05)  {\includegraphics[width=.51\textwidth]{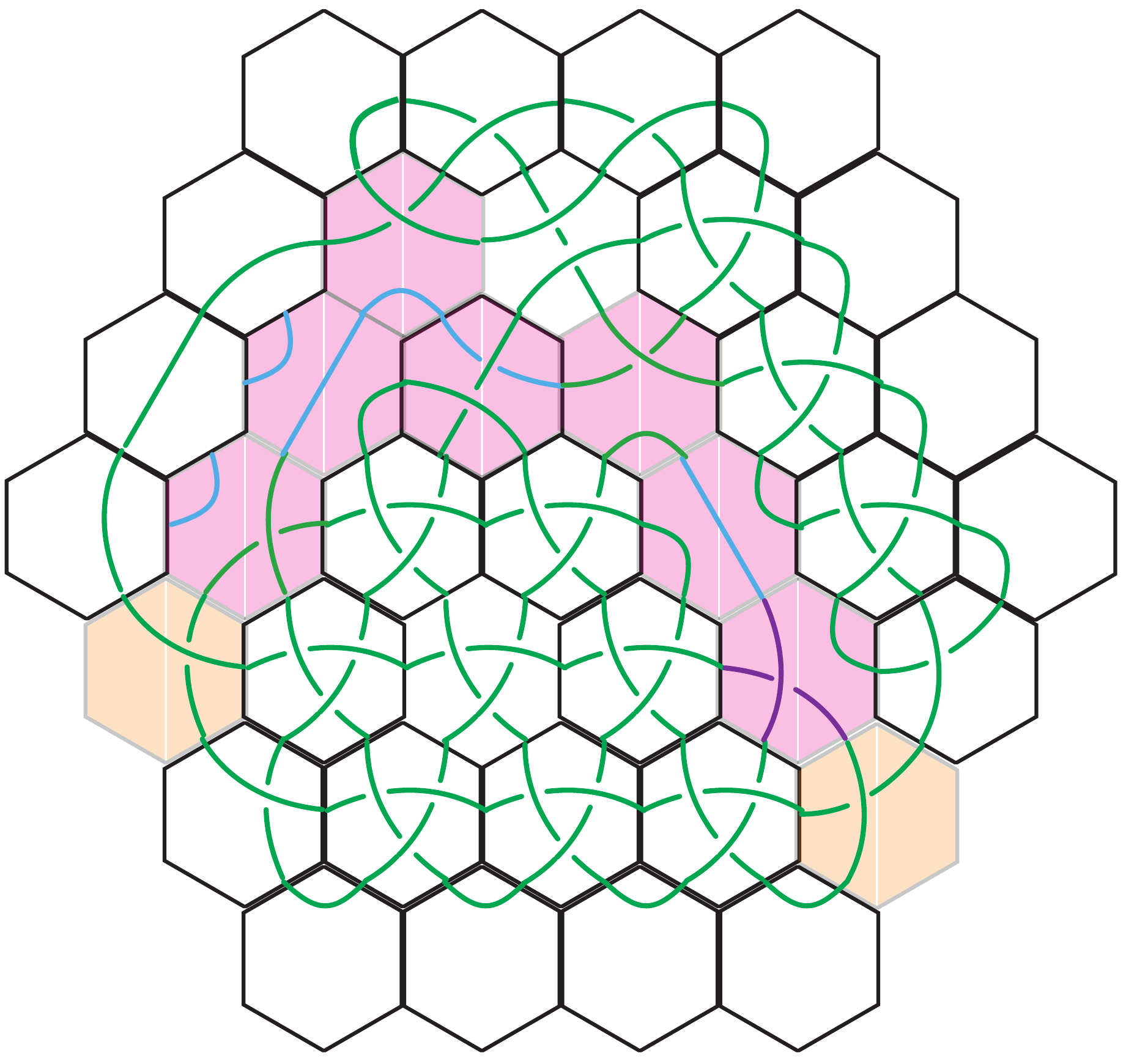}}
	\put(.6,4.1){$K^2$}
	\put(5.9,4.1){$K^3$}
	\put(8.75,1.66){$T_{5,5}$}

	\end{picture}
	\caption{$L^2$ in Figure~\ref{fig:LA2} had multiple components.  On the left $K^2$ bands components of $L^2$ together to form a knot. On the right we form a final knot $K^3$ from $K^2$ making sure that $K^3$ has at least as many crossings as $K^1$ did, but $K^3$ has one less arc in its complement than $K^1$ did. }
	\label{fig:KA2}
	\label{fig:LA3}
	\label{fig:KA3}
	\end{figure}
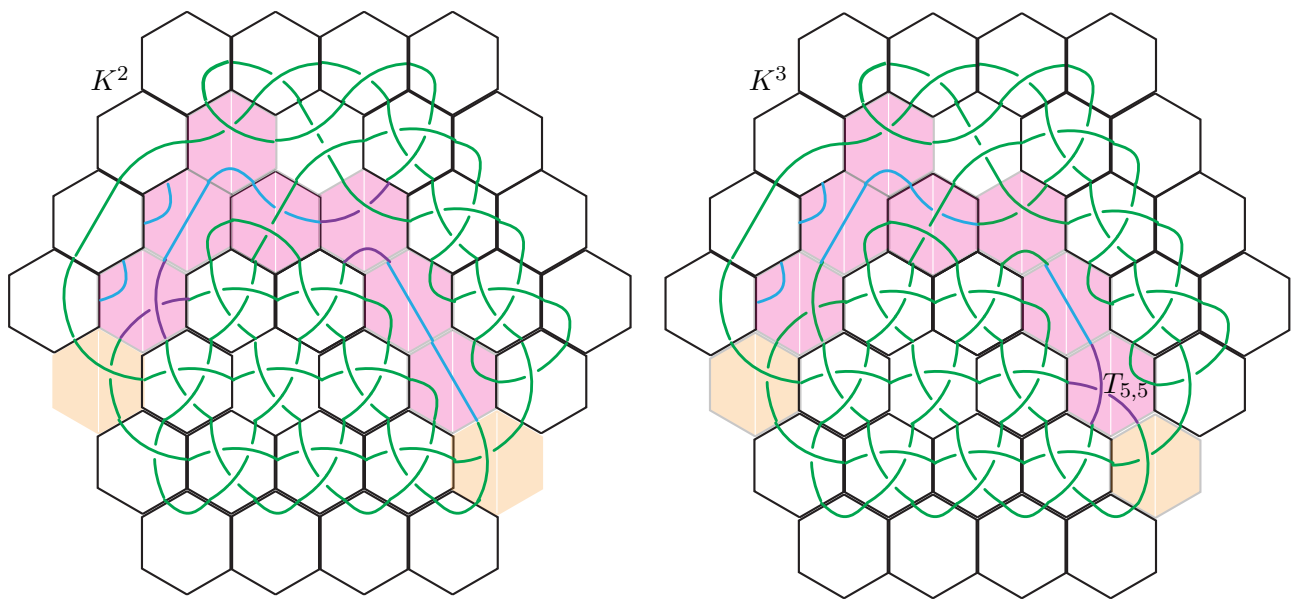

	\begin{prop}
		The knot projection with the largest number of crossings that can be obtained on an $r$-mosaic can be chosen so it has no arcs or loops in its complement.  This is true for  rectangular $r$-mosaics as well as standard, semi-enhanced, and enhanced hexagonal $r$-mosaics.
		\label{prop:cross}
	\end{prop}

	The process we develop here  reproves  Corollary 7.16  in 
	\cite{hk} for rectangular mosaics vastly reducing the length of the proof of the main result in that paper and also leads to a brand new result for standard, semi-enhanced, and enhanced hexagonal mosaics.
	
	\medskip

	\emph{Proof of Proposition~\ref{prop:cross}:}
	In Section~\ref{sec:fam}  we formed a knot called $A_r$ by smoothing  crossings of a saturated link $L_r$ for all rectangular and hexagonal settings. In each case we want to show that it is impossible to find a knot on an $r$-mosaic with higher crossing number than $A_r$.
		We show that for the various types of hexagonal $r$-mosaics with the single exception of $\ddddot{A_3}$ and for rectangular $r$-mosaics with even $r$ there can never be a knot on an $r$-mosaic with more crossings than $A_r$ (regardless of crossing number).  
		
		For hexagonal mosaics, the aberrant and less interesting case is the standard  hexagonal 3-mosaic where  it is possible to find a knot isotopic to $\ddddot{A_3}$ 
 with one nugatory crossing.  This is not true for $\widehat{A_3}$ or $\hat{A_3}$.  For rectangular $r$-mosaics we are skipping the case of odd $r$ since it is easy and covered in \cite{hk}, but note that we can embed a knot isotopic to $\overline{A_r}$, but with up to two nugatory crossings, so in those cases there will be embeddings of knots with more crossings than $A_r$, but they are always are reducible, isotopic to $A_r$ and there never exist knots with higher crossing number than $A_r$ in any of these settings.
	The construction below is nearly identical for rectangular, standard, and semi-enhanced hexagonal mosaics and we will  point out the few minor differences as we go along.  There is extra work for enhanced hexagonal mosaics as seen below.  
	

	By definition the complement of a link mosaic consists
	of loops and arcs.   Let $s \geq 0$ be the number of loops in the complement and $w \geq 0$ be the number of arcs for a given mosaic. 
	Of all knots one can embed on an $r$-mosaic we will choose $K^1$ from the collection of embedded knots that maximize the number of crossings (not the crossing number, but the number of crossings in the specific embedding). If there are no knots with more crossings than
	$A_r$ we are done, so further assume the choices for $K^1$ have more crossings than $A_r$.  
	 Over all options for $K^1$, examine the ones whose complements minimize $s$ and then among those possible choices pick one that then minimizes $w$ (in other words choose $K^1$ and its complement so that the number of crossings for the knot is maximized, but within those options so that the ordered pair $(s,w)$ is minimized lexicographically).  We will first show that this implies $s=0$ and then that $w=0$, too.  The argument for $s=0$ is the same for all forms of hexagonal mosaics and essentially the same for rectangular mosaics.  The argument that $w=0$ will require adjustments going from one case to the other.

	\medskip
\noindent	\emph{Showing that $s=0$:}
	\medskip

\begin{figure}[tpb]
		\centering
		\includegraphics[width=.36\textwidth]{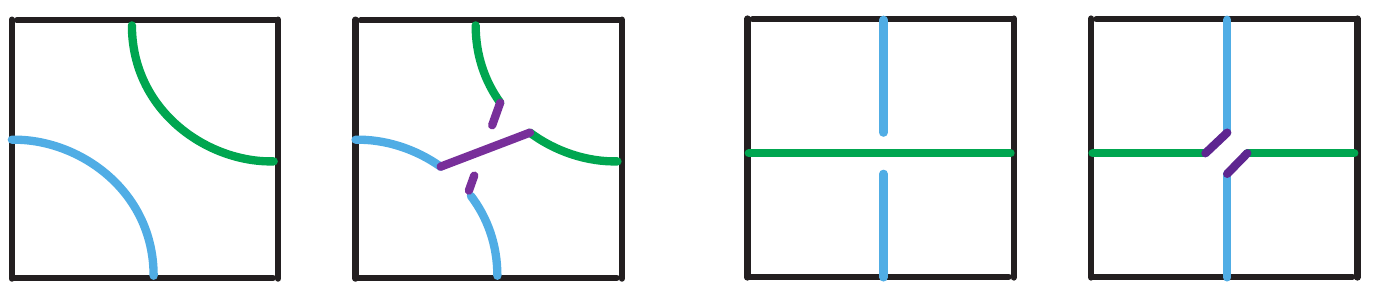}
		\includegraphics[width=.91\textwidth]{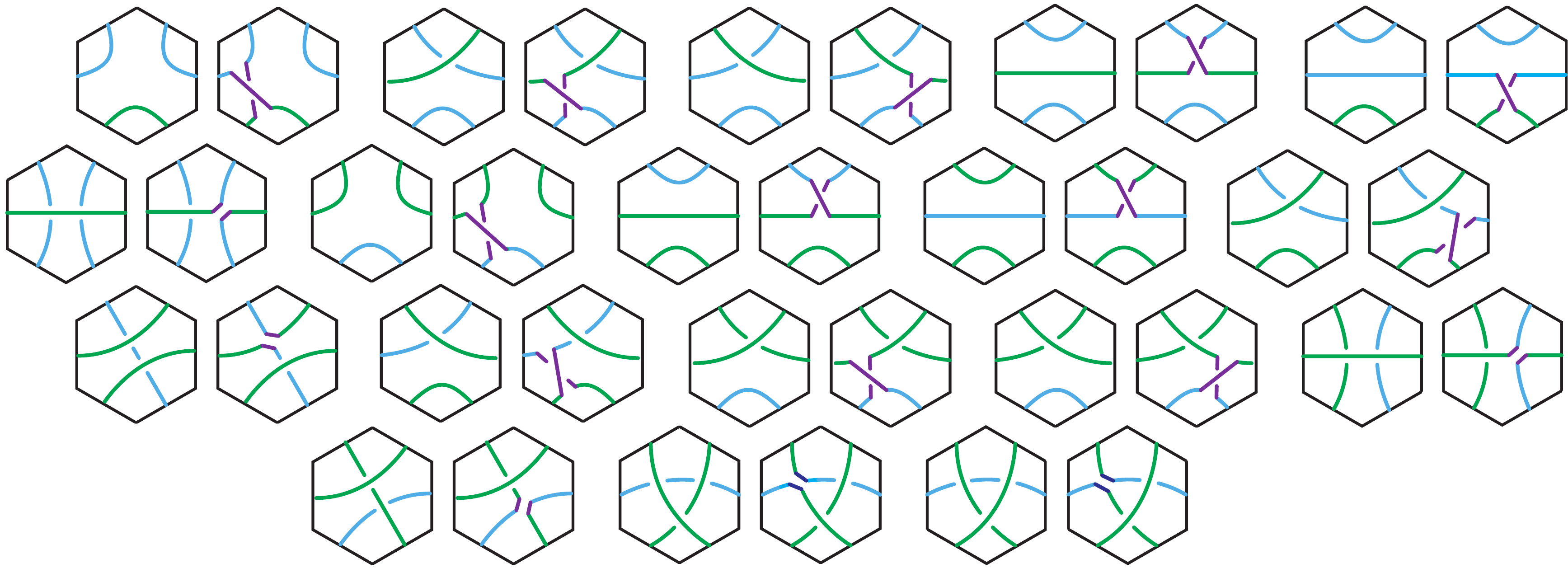}
		\caption{We can combine an arc of the complement with an arc of the link (or another arc of the complement) by banding them together or smoothing a crossing between the knot and the complement as shown above for the complement with the link. The top row shows how to do this for rectangular mosaics and the bottom rows show it for hexagonal mosaics.  Each time the knot is drawn in green and the complement in blue and to the right we see the arcs banded together with purple arcs.  Note that if the two arcs being combined do not already cross each other we use the band to gain an extra crossing.}
		\label{fig:comb}
		\label{fig:comb2}
	\end{figure}

	If $s>0$ (there are loops in the complement), then we examine one of those loops which we call $m_1$. We will show that you can take the connect sum of $m_1$ with $K^1$ or with another component of the complement leading to a contradiction.  We know $m_1$ must, of course, either enter a tile with another component of the complement or a tile containing a portion of $K^1$ (or possibly both).  First imagine it intersects a tile containing an arc of $K^1$.   We then can band together $m_1$ with $K^1$ to get a new knot with one less loop in its complement.  This tile replacement process is shown in Figure~\ref{fig:comb}. In the figure a tile of $K^1$ drawn in green and the complement drawn in blue is pictured and directly to its right is a way to combine the two together.   The result is the connect sum of $K^1$ and the (unknotted) loop from the complement. The connect sum is a new mosaic and the knot it represents is isotopic to $K^1$.
	The new knot mosaic still contains all the crossings of $K^1$ (and often even more crossings).  This shows there was a knot mosaic with at least as many crossings as $K^1$ and fewer loops in its complement contradicting our earlier minimality assumption.

\begin{figure}[tpb]  \centering
	{\includegraphics[width=.8\textwidth]{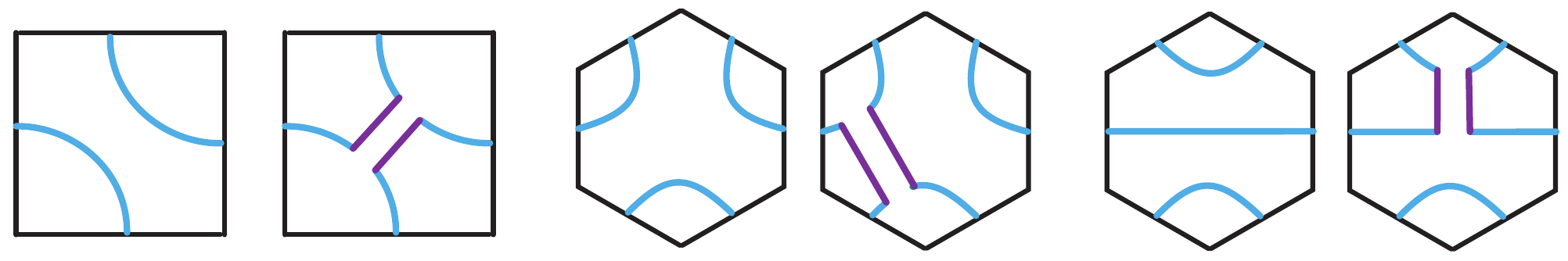}}
		\caption{We can combine two arcs of the complement when the knot itself is disjoint from the tile.  If the component of the complement, $m_1$, never enters a tile containing the knot there must be a tile that it enters which also contains a different component of the complement, $m_2$, and we band $m_1$ to $m_2$.   On the left we see the only possible square tiles with only arcs of the complement and a way to connect them and on the right we see the two possibilities for hexagonal tiles.}
	\label{fig:compcomb}
\end{figure}

	We established above that if $s>1$ then there is a loop $m_1$ of the complement, and that $m_1$ cannot enter any of the tiles containing $K^1$ without a contradiction.  Thus it must enter a tile containing another component of the complement (one can see this, for example, because a tile containing a global minimum of $m_1$ must also contain an arc not from $m_1$ below the minimum, either from $K^1$ or another component of the complement and since the tile is disjoint from $K^1$ the second arc must be from another component of the complement).  Call the second component of the complement $m_2$, but note that $m_2$ may be a loop, but also might be an arc of the complement.   For hexagonal $r$-mosaics we see one of the first two tiles in top row of Figure~\ref{fig:comp}  and for rectangular $r$-mosaics we see the tile labeled 1 in Figure~\ref{fig:squarecomp}
	since those are the only mosaic tiles that are disjoint from the knot and  we assumed the knot does not enter the tile.  If  the arcs of $m_1$ and $m_2$ are separated by an arc of a third component of the complement $m_3$  then we retake $m_2$ so the arcs are adjacent in that tile.  Now we band $m_1$ and $m_2$ together as in Figure~\ref{fig:compcomb}, taking a connect sum of either two loops or a loop and an arc.  This yields a different complement for the same mosaic (as we pointed out earlier the complement is not well defined).  The knot is unchanged, but the number of loops in the complement dropped by one, so we get our final contradiction of $s$ being greater than 0.  Thus the complement contains no loops.

		\medskip
	\noindent	\emph{Showing that $w=0$:}
	\medskip

	We now argue that $w=0$, too, showing the complement can have no arcs either. For the sake of contradiction, assume $w>0$.  Recall $K^1$ has as many crossings as any knot embedded on an $r$-mosaic but as few arcs as possible in its complement (but as we showed above no loops in the complement).   
	We will use an outermost arc of the complement (we define the outside of an arc of the complement and outermost arcs of the complement in the next paragraph) to build a sequence of links and knots until we end with a knot $K^3$ which has at least as many crossings as $K^1$ and fewer arcs in its complement, a contradiction of the assumption that $w>0$.   
	
		As always the arcs of the complement lie on the interior tiles of the mosaic with both endpoints of the boundary corona like the blue arc of the complement $a$ on the right in Figure~\ref{fig:grid}.   Let the set of interior tiles containing $a$ be called $A$.  Let the two tiles of the boundary corona which contain the endpoints of $a$ be called $\tilde{A}$  (if the endpoints of $a$ lie on the same tile then $\tilde{A}$ is a single tile, of course).  $\tilde{A}$ divides the boundary corona into two pieces: we allow one of the pieces to be trivial if $\tilde{A}$ consists of a single tile or adjacent tiles.   If one side of $\tilde{A}$  in the boundary corona contains fewer boundary tiles than the other does (or no tiles at all), the smaller side is called the {\em outside} of $a$.  If $a$ is a diameter in the sense that the number of boundary tiles on both sides is equal then either side may be called the outside.  If there are acrs in the complement, there must be some arc $a$ with no arcs of the complement outside of it and we will call this {\em an outermost arc of the complement}.
	
 	The tiles of our mosaic break down into 6 groups, the first two sets are   $A$ and $\tilde{A}$.   In the mosaic on the right in Figure~\ref{fig:grid} $A$, of course, is the set of tiles containing the blue arc labeled $a$ and are shaded pink. The tiles of $\tilde{A}$ consists of the two tiles $T_{5,1}$ and $T_{6,5}$ labelled and shaded orange.  The interior tiles not containing $a$, but outside $a$ are called $O$.  In the same figure these are the interior tiles below $A$. The interior tiles not containing $a$ but inside $a$ are called $I$ (all the interior tiles above $A$ in Figure~\ref{fig:grid}).

	The portion of the boundary corona outside $A \cup \tilde{A}$ is called $\tilde{O}$ (in Figure~\ref{fig:grid} $\tilde{O}$ consists of tiles $\{T_{6,1}, T_{7,1}, T_{7,2}, T_{7,3}, T_{7,4} \}$), and the boundary tiles inside $A \cup \tilde{A}$ are called $\tilde{I}$. Our algorithm will produce mosaics $K^1, L^1, L^2, K^2$, and $K^3$. All of these mosaics  will perfectly agree on the tiles $I$ and $\tilde{I}$, but we will change some of the tiles in $A$, $\tilde{A}$, $O$, and $\tilde{O}$.

	Now that the outside of an arc is defined, in our proof we choose $a$ to be an outermost arc of the complement of $K^1$.

	\medskip
	\noindent	\emph{Step 1: form $L^1$ by adding $a$ to $K^1$:}
	\medskip

To start the construction let $J$ be the set of interior arcs of $K^1$ ($K^1$ intersected with the interior tiles of the mosaic)   
	and form $L^1$, a link whose interior arcs consist of $J \cup a$.  Once the interior arcs are defined we only have a few ways to form a link through the boundary tiles: exactly two in the case of rectangular mosaics and standard hexagonal mosaics but a few more in the enhanced and semi-enhanced cases.  We can always choose to connect so that either $K^1 \cap \tilde{I}$ or $K^1 \cap \tilde{O}$ remains unchanged as we form $L^1$, but the other side will be altered.  We choose $L^1$ to be the link formed by connecting through the boundary tiles in the way that leaves $K^1 \cap \tilde{I}$ unchanged, but alters $K^1 \cap \tilde{O}$ and $K^1 \cap \tilde{A}$.    We see this step as we go from the knot $K^1$ in Figure~\ref{fig:KA} to the link $L^1$ depicted on the left in Figure~\ref{fig:LA1}.  In the rectangular case and standard hexagonal case
	this choice completely determines $L^1$, but in the enhanced and semi-enhanced cases we have some choice. 
	In the semi-enhanced case make the exact same choices on $\tilde{A} \cup \tilde{O}$ as we did for the standard case so the two links would look identical on $\tilde{A}  \cup \tilde{O}$.
	 In the enhanced case we choose to fill in $\tilde{A}$ and $\tilde{O}$ in the unique manner that has as many crossing tiles as possible on $\tilde{A} \cup \tilde{O}$.

	In each of the settings we formed a new link from $K^1$ by adding $a$ to it, but in the enhanced case we might lose some crossings in the boundary tiles and if so we will have to make sure losses are compensated for in another part of the construction.
	 We see crossings vanish in boundary tiles  $T_{7,2}$ and $T_{7,3}$ when we form $L^1$ in Figure~\ref{fig:LA1} from $K^1$ in Figure~\ref{fig:grid} (in that specific case we also gain crossings in boundary tiles $T_{5,1}$, $T_{6,1}$ and $T_{6,5}$,    more than compensating for the lost crossings, but this need not always be the case).	   In the rectangular case, the standard hexagonal case, and the semi-enhanced case all the original crossings of $K^1$ will persist in each of the stages of the construction of $K^3$ showing $K^3$ clearly has at least as many crossings as $K^1$.  The enhanced hexagonal case requires a little more attention.
	The entire construction is shown for enhanced hexagonal mosaics in Figures~\ref{fig:grid}, \ref{fig:LA1}, and \ref{fig:LA3}.

\medskip
	\noindent	\emph{Step 2: form $L^2$ from $L^1$ by saturating the tiles of $L^1 \cap O$:}
	\medskip

	Since $a$ is outermost and there are no loops in the complement, the tiles of $O$ are disjoint from the complement and thus $K^1$ and now $L^1$ used all of the connection points in $O$.  We can thus freely replace these tiles of $L^1$ with saturated tiles: tile 26 or 27 in the hexagonal setting and the crossing tile in the rectangular case.  Call this new link  $L^2$.  This step is depicted in Figure~\ref{fig:LA1} as we form $L^2$ on the right from $L^1$ on the left by adding  three crossings on tile $T_{5,3}$.    We have now ensured the tiles of $L^2 \cap (O \cup \tilde{O})$ are fully saturated having as many crossings as possible.
	Clearly $L^2$ has at least as many crossings as $L^1$ and more unless $L^1$ was already saturated on these tiles making this step trivial.  Because we are only currently concerned with the number of crossings in the mosaic, crossing changes are unimportant so we may choose the assignments of over and under at each crossing to ensure that $L^2 \cap O \cup \tilde{O}$ is identical to  $L_r \cap O \cup \tilde{O}$ using the same type of argument as we did in the proof of Lemma~\ref{lemma:iso}.

\medskip
	\noindent	\emph{Step 3: form a knot $K^2$ from the link $L^2$:}
	\medskip

	We now have moved $a$  out of the complement and have saturated the tiles of $O$, but we may have formed a link with multiple components and need to get back to a knot.
	If $L^2$ has only one component we rename it $K^2$ and proceed to the next step.
	
	  If $L^2$ has multiple components call the one containing $a$, $C_a$ and choose a component not containing $a$ and call it $C_b$. Since the tiles of $K^1 \cap (I \cup \tilde{I}) = L^2 \cap (I \cup \tilde{I})$ and $K^1$ is a single knot, $C_b$ cannot strictly contained in $I \cup \tilde{I}$.  Since 
	$L^2 \cap (O \cup \tilde{O}) = L_r \cap (O \cup \tilde{O})$ and
	$L_r$ never has an entire component outside a properly embedded arc on the interior tiles of a hexagonal mosaic, $C_b$ must intersect the tiles $A \cup \tilde{A}$ at least twice in the hexagonal case (in general you can draw a properly embedded arc on $L_r$ that is disjoint from one or more of the components of $L_r$, but those missed components will lie on the inside of the arc, not the outside).  
	
	Here is the only time in the paper where we have to take special care in the rectangular case that is not also necessary in the hexagonal cases.  In the even rectangular case $\overline{L_r}$ could have an entire component outside of an arc $a$ if $a$ is long enough to hit opposite sides of the mosaic, but that never happens in our situation since $\overline{A_r}$ only requires smoothing $r-3$ crossings of $\overline{L_r}$ so  $a$ hits at most $r-3$ tiles and thus must have endpoints either on the same side of the boundary of the interior of the mosaic or adjacent sides eliminating this risk.  In the hexagonal setting we do not need to use the bound on the length of $a$ to ensure there is no component strictly contained on $O \cup \tilde{O}$ or $I \cup \tilde{I}$ since this cannot happen outside $a$ no matter how long $a$ is.

	As we know $a \subset A$ and has its endpoints on $\tilde{A}$: let $a'$ be $a$ extended into $A \cup \tilde{A}$ so that $a \subset a' \subset C_a$.
	If $C_b$ crosses $a'$ then it does so twice and these are new crossings that did not exist in $K^1$ and we can smooth one of them netting one extra crossing (the one which we did not smooth) while decreasing the number of components by one.  We do this in Tile $T_{3,4}$ in $L^2$ on the right in Figure~\ref{fig:LA1}  as one step in forming $K^2$, dropping from 3 to 2 components in Figure~\ref{fig:LA3}.

	If $C_b$ does not cross $a'$ then it must enter $A \cup \tilde{A}$ and turn back instead of crossing $a'$. We now argue that even if $C_b$ enters $\tilde{A}$ and turns back it must also enter $A$ and turn back.  This is certainly true vacuously in the rectangular case and the enhanced case because in the rectangular case any arc that enters $\tilde{A}$ must connect to $a$ and thus must be part of $C_a$ not $C_b$.  In the enhanced case the construction ensured that such an arc must either connect to $a$ or cross $a'$ so again it cannot turn back in $\tilde{A}$.  In the standard hexagonal case and the semi-enhanced case it is possible for an arc to enter $ \tilde{A}$ and turn back, not connecting to $a$ or crossing $C_a$, but because in both cases our construction ensures that $L^2$ looks like $\ddddot{L_r}$ on $O \cup \tilde{O}$ and since every component of $\ddddot{L_r}$ that intersects $O \cup \tilde{O}$ also intersects $A$ at least twice we know $C_b$ must enter $A$ at least twice.  Since earlier we took care of the case where $C_b$ crosses $a$ we may conclude that $C_b$ enters $A$ at least twice and turns back without crossing $a$ either time.   In that case we band the new component to $a$ adding in a crossing as in  Figure~\ref{fig:comb2} reducing our total number of components by one and increasing our crossings.  Repeat until we have a single knot and call this knot $K^2$.

   We did this to the tile $T_{4,2}$ in $L^2$ on the right in Figure~\ref{fig:LA1} as the other step in forming $K^2$ in Figure~\ref{fig:LA3}. In Figure~\ref{fig:comb2} the blue arcs represent $a \subset C_a$ passing through the tile and the green arcs represent an arc of $C_b$.  We band sum together $C_b$ and $C_a$ within that tile as shown in Figure~\ref{fig:comb2} to form a new link with one less component than $L^2$ had, but with one  crossing that was not part of  $K^1$.    If $L^2$ had $t+1$ components, $K^2$ has at least $t$ more crossings on the interior tiles than $K^1$ did.  In the standard mosaic case or the rectangular case we have found a knot with more crossings than $K^1$ since $L^2$ had more crossings than $K^1$, a contradiction.

%

	  \medskip
	  \noindent	\emph{Step 4: Finally, form a final knot $K^3$ out of $K^2$:}
	  \medskip

	As we mentioned above we are working to form a new knot $K^3$ from $K^1$ and the first steps were to form $L^1$ then $L^2$ then $K^2$.  Clearly none of the stages ever have fewer crossings on the interior tiles than $K^1$ and all of them have fewer arcs in the complement.  In the rectangular, standard hexagonal, and semi-enhanced hexagonal settings there are no crossings in boundary tiles to worry about so we can let $K^3=K^2$ and we have our contradiction to minimality so from here on we restrict to the enhanced hexagonal setting.  
	
	 In the enhanced hexagonal setting a worry is we may have lost crossings in the boundary tiles forming $L^1$ (which we will now call $\widehat{L^1}$ since we are only worried about enhanced hexagonal mosaics).  To make sure this is not a problem we examine the boundary tiles of $\widehat{L^1}$.  Since $\tilde{I}$ is unchanged losses can only happen on $\tilde{A}$ or $\tilde{O}$.  $\tilde{A}$ originally had the end points of $a$ on it, so $\widehat{K^1}$ could not have had crossings on those tiles, so we only have to worry about $\tilde{O}$.  It is, of course, possible that $\tilde{O}$ has more crossings than it did before, but this would only help us so we look at how many crossings we could possibly lose.  
	 
	    \begin{lemma}
	   Boundary crossing tiles cannot occur on adjacent sides of the boundary hexagon of a mosaic outside of the an outermost arc such as $a$.
	  \label{lemma:boundarysat}
	\end{lemma}

\begin{proof}
This follows directly from Lemma~\ref{lemma:boundarycrossings} and the definition of an outermost arc of the complement.
\end{proof}

	\begin{lemma}
	The total number of crossings of $\widehat{K^1} \cap \tilde{O}$ can exceed that of $\widehat{L^1} \cap \tilde{O}$ (which equals the number for $\widehat{K^2} \cap \tilde{O}$) by no more than $r-2$.
	
	\end{lemma}

	\begin{proof}
	Obviously $a$ has to 	
	be relatively long to be problematic or    $\widehat{K^1} \cap \tilde{O}$ could not have $r-2$ crossings on it.  
	  There are two critical points for the number of possible lost crossings on $\tilde{O}$ as we go from $\widehat{K^1}$ to $\widehat{L^1}$
	  which form worst case scenarios. The first is when $a$ runs from one corner tile to an adjacent corner tile with each of the tiles of $\widehat{K^1} \cap \tilde{O}$ containing crossings.  In this case we lose those $r-2$ crossings when we form $\widehat{L^1}$.  If $a$ were slightly longer and cut off part of a second or third side $\widehat{K^1} \cap \tilde{O}$, then $\tilde{O}$ would have to contain boundary tiles without crossings   and these become crossing tiles of $\widehat{L^1} \cap \tilde{O}$ decreasing the total number of crossings lost.  The other maximum occurs when $a$ is a diameter running from one corner of the mosaic to another corner and $\widehat{K^1} \cap \tilde{O}$ consists of two sides each containing $r-2$ crossings and a third without crossings.  In this case as we form $\widehat{L^1}$, on $\widehat{L^1} \cap \tilde{O}$ we lose the $2(r-2)$ crossings, but the side that started without crossings gives us $r-2$ new crossings for a net loss of $r-2$ crossings total.

	  \end{proof}
	
	We now observe that by construction $\widehat{K^2}$ is saturated on $O \cup \tilde{O}$ and $\widehat{K^2} \cap (O \cup \tilde{O})=\widehat{L^2} \cap (O \cup \tilde{O})$ has been chosen to exactly match $L_r \cap (O \cup \tilde{O})$. 
	$\widehat{K^2} \cap (O \cup \tilde{O})$ can be thought of as a series of arcs with endpoints on $A \cup \tilde{A}$.  Within this series there may be a few arcs each with one  or both endpoints on $\tilde{A}$, but all the rest of the arcs have both endpoints on $A$.  We will show that if there are $j$ arcs with both endpoints on $A$ and $n$ crossings lost in $\tilde{O}$ then $j \geq n$.  We will then use these arcs to cancel out any lost crossings.

	Let $n$ equal
	the  number of crossings of $\widehat{K^1} \cap \tilde{O}$ minus the number of crossings on $\widehat{K^2} \cap \tilde{O}$.  We showed above $n \leq r-2$.   Let    $\{ e_1, e_2 \dots e_j \}$ be the arcs in $\widehat{K^2} \cap (O \cup \tilde{O})$ where each $e_i$ has both endpoints on $A$.  We will argue  $j \geq n$.  Each $e_i$ is used to create a new crossing that wasn't contained in $\widehat{K^1}$ compensating for any crossings lost in the boundary tiles.  This breaks down into cases based on how many edges of the boundary hexagon intersect  $\tilde{A} \cup \tilde{O}$.

	If $\tilde{A} \cup \tilde{O} $ is contained on a single edge then all the tiles of $\widehat{K^2} \cap \tilde{O}$  are either crossing tiles or none are.  If they are all crossing tiles then $n$ is negative and we are done so we assume that none of them are.  In that case each of these tiles intersects an edge of $\widehat{K^2} \cap (O \cup \tilde{O})$ which hits $A$ twice yielding $\{ e_1, e_2 \dots e_j \}$.	 

If $\tilde{A} \cup \tilde{O} $ intersects two sides of the boundary hexagon then one of those sides contains crossing tiles and the other contains only non-crossing tiles by Lemma~\ref{lemma:boundarysat}.  Again each of the lost crossings in $\tilde{O}$ came from the new non-crossing tiles, and again each of these tiles intersect an edge of $\widehat{K^2} \cap (O \cup \tilde{O})$ which hits $A$ twice.  If we rotate the mosaic so that the side now containing the non-crossing tiles is at the bottom then each $e_i$ will contain a local minimum on one of these tiles and will look like an edge of $\widehat{L_r} \cap (O \cup \tilde{O})$ with $\widehat{L_r}$ rotated so that its bottom edge contains no crossing tiles.

If $\tilde{A} \cup \tilde{O} $ intersects three sides of the boundary hexagon then there will be $r-2$ edge tiles in $\tilde{O}$ running between two corner tiles in $ \tilde{O} $ regardless of whether these are crossing tiles or not they intersect $r-2$ arcs of $\widehat{K^2} \cap (O \cup \tilde{O})$ with both endpoints on $A$ giving $e_1, e_2 \dots e_{r-2}$ and since we showed  $n \leq r-2 \leq j$  we have the desired set of arcs.  The other sides may well contribute more such arcs, but that is just a bonus and not needed for the argument.

Finally it is possible that $\tilde{A} \cup \tilde{O} $ intersects four sides of the boundary hexagon.  In that case there are two entire boundary edges contained in $\tilde{A} \cup \tilde{O}$, one of which contains $r-2$ crossing tiles and the other of which contains none.  If we examine the boundary edge with no crossing tiles we see at least $r-3$ arcs of $\widehat{K^2} \cap (O \cup \tilde{O})$ with both endpoints on $A$.  It is possible that (at most) one of the tiles yields an arc with an endpoint on $\tilde{A}$ or we would already have the desired $r-2$.  One can easily argue that the side with crossings on it also has multiple arcs of  $\widehat{K^2} \cap (O \cup \tilde{O})$ with both endpoints on $A$ pushing us well past the desired bound of $r-2$ edges in our set $\{ e_1, e_2 \dots e_j \}$.

\begin{figure}[tpb]  \centering
	\setlength{\unitlength}{0.1\textwidth}
	\begin{picture}(10,1)
	\put(1.95,0.22){\includegraphics[width=.61\textwidth]{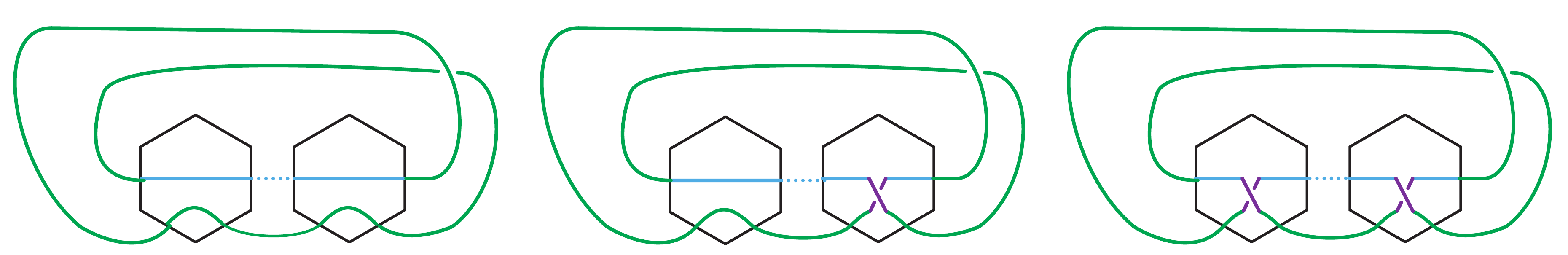}}
		\put(5.0,.1){$e_i'$}
	\put(5.0,.55){$a$}

		\put(2.95,.1){$e_i'$}
	\put(2.95,.55){$a$}

		\put(7.05,.1){$e_i'$}
	\put(7.05,.55){$a$}
		\end{picture}

	\caption{Here on the left we see an edge $e_i'$ entering the tiles of $A$, but not crossing $a$ (depicted in blue).  Adding a crossing to either tile individually as we do  in the middle picture for the tile on the left turns the knot into a link of two components, but in that case adding a crossing to both tiles yields a knot with two more crossings than we started with as we see on the far right. }
	\label{fig:twocross}
\end{figure}
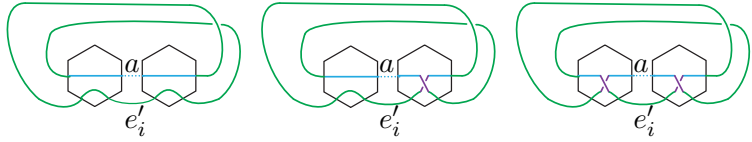

Now we have shown that in each case $j \geq n$ and we have at least as many arcs of $\widehat{K^2} \cap (O \cup \tilde{O})$ with both endpoints on $A$ as crossings we lost in the boundary forming $\widehat{K^2}$.  Examine these arcs one at a time starting with any given $e_i$.  We already gained a crossing back in step 3 for each $e_i$ that was part of a component of $\widehat{L^2}$ not containing $a$ when it was banded together  to form $\widehat{K^2}$ so we need only find a way to add crossings when $e_i$ is one of the other edges.   By definition $e_i \subset \widehat{K^2} \cap (O \cup \tilde{O})$. Let $e_i'$ be the arc of $\widehat{K^2} \cap   (O \cup \tilde{O} \cup A)$ containing $e_i$ and extending it into $A$ by one tile one each end.  
 If $e_i'$ crosses $a$ at least once after entering $A$ then this netted a new crossing in $\widehat{K^2}$ that didn't exist in $\widehat{K^1}$ and that compensates for the lost crossing in the boundary.
 
 Thus our only remaining concern is the case where $e_i$ enters $A$ at both endpoints but $e_i'$ does not cross $a$.  We could then use one of the band sum replacements shown in Figure~\ref{fig:comb} to  connect $e_i$ and $a$.  We use such a move on $\widehat{K^2}$ in Figure~\ref{fig:LA3} in tile $T_{5,5}$ to go from $\widehat{K^2}$ on the left to $\widehat{K^3}$ on the right.  Since in general $e_i'$ enters $A$ twice we have two choices where we can do this.  If banding either one turns $\widehat{K^2}$  into another knot as it did in the Figure~\ref{fig:LA3} then we make this choice getting one extra crossing again negating the lost crossing in the boundary tile.  It is possible that both choices for bands individually turn $\widehat{K^2}$ into a link of two components.  This can only be true if the arcs of $\widehat{K^2} - (a \cup e_i')$ cross an odd number of times as they do in Figure~\ref{fig:twocross}.  
 In that case, however, as we see in the figure adding in both possible crossing bands like we do on the far right yields a knot because the first band creates a link and the second band joins the two link components back together.   We follow the process above for each edge in $\{ e_1, e_2 \dots e_j \}$ and call the result $\widehat{K^3}$.  $\widehat{K^3}$ has at least as many crossings as $\widehat{K^1}$ and one less arc in the complement contradicting minimality.

 Thus we have shown that $w=0$ and we may thus assume that the complement of $K^1$ is trivial.

 \[
\pushQED{\qed} 
 \qedhere
\popQED
\]

	Proposition~\ref{prop:cross} in turn leads directly to the following Corollary.

	\begin{cor}
		For rectangular $r$-mosaics and for standard, semi-enhanced, and enhanced hexagonal $r$-mosaics for $r \geq 2 $,
		the knot projection with the largest number of crossings that can be obtained on an $r$-mosaic can be obtained by smoothing  the crossings of a saturated mosaic.
		\label{cor:smooth}
	\end{cor}

	\begin{proof}
		Since we may assume there are no arcs or loops in the complement we see immediately that the maximal crossing knot may be obtained  by smoothings from a mosaic board whose interior is saturated. Note that the corollary is trivially true for hexagonal $2$-mosaics since   $\ddddot{L_2}= \hat{L_2}= \widehat{L_2}= \ddddot{A_2}= \hat{A_2}= \widehat{A_2}$.
		For an enhanced mosaic if the interior tiles are saturated all corner tiles must be Tile 2 from Figure~\ref{fig:tiles}.  Three of the boundary edges cannot have crossings on them and must also be made of Tile 2 by Lemma~\ref{lemma:boundarysat}.  The remaining 3 edges must be made up of Tiles 10 or 11 if they contain  a crossing and Tiles 5 or 6 if not.  Note that Tiles 5 and 6 can be obtained from Tiles 11 or 12 by smoothing.  Thus if the complement is trivial we already knew that the interior tiles could be obtained by smoothing three crossing tiles, but now we know any boundary tiles which could have contained a crossing, but do not can also be thought of as smoothings of crossing tiles so we may in every setting assume we started with a saturated link and smoothed crossings.		
	\end{proof}

	We now come to the main results of the paper. We first reprove the following theorem from \cite{hk} more efficiently.
	
	\bigskip
	
	\begin{thm}
		 
	\noindent (Theorem 8.2 in \cite{hk}) 
	Given a rectangular $r$-mosaic with $r>3$ and any knot $K$ that is projected onto the mosaic, the crossing number $c$ of $K$ is bounded above by the following:\\
	\begin{equation*}
	c \leq
	\begin{cases}
	(r - 2)^{2} - 2 & \quad \text{if $r = 2k + 1$}\\
	(r - 2)^{2} - (r - 3) & \quad \text{if $r = 2k$.}
	\end{cases}
	\end{equation*}

	Note that $\overline{A_r}$ achieves this crossing number for $r > 3$, proving the bound is sharp.  

	\label{thm:rect}
\end{thm}
	
	\bigskip
		In the rectangular theorem above we start with $r>3$ because one cannot find a non-trivial knot on a rectangular $r$-mosaic with $r \leq 3$ so the result is trivially 0 when $r \leq 3$).  We now state the new result for hexagonal $r$-mosaics and will prove both theorems together below.
		
		\bigskip
		
	\begin{theorem}

		Given a standard hexagonal $r$-mosaic and any knot $K$ that is projected onto the mosaic, the crossing number $c$ of $K$ is bounded above by the following:\\

		\begin{equation*}
		c \leq
		\begin{cases}
		3=9r^2-28r+23 & \quad \text{if $r = 2$}\\
		 19= 9r^2-28r+22 & \quad \text{if $r = 3$}\\
		9r^2-28r+23 & \quad \text{if   $r >3$.}
		\end{cases}
		\end{equation*}
\\

		Given a semi-enhanced hexagonal $r$-mosaic and any knot $K$ that is projected onto the mosaic, the crossing number $c$ of $K$ is bounded above by the following:\\

		\begin{equation*}
		c \leq
		\begin{cases}
		3  & \quad \text{if $r = 2$}\\
		9r^2-27r+21 - (\lceil \frac{r}{2} \rceil - 1) = 9r^2-27r+22 - \lceil \frac{r}{2} \rceil& \quad \text{if   $r >2$.}
		\end{cases}
		\end{equation*}
\\

Given an enhanced hexagonal $r$-mosaic and any knot $K$ that is projected onto the mosaic, the crossing number $c$ of $K$ is bounded above by the following:\\

		\begin{equation*}
		c \leq
		\begin{cases}
		3 & \quad \text{if $r = 2$}\\
		9r^2-25r+15 & \quad \text{if $r >2$.
		}
		\end{cases}
		\end{equation*}
	\\

		Note that $\ddddot{A_r}$ $\hat{A_r}$, and $\widehat{A_r}$ achieve the exact bounds for $r \geq 2$, proving the bounds are sharp.
		\label{maxcross}
		\label{thm:maxcross}
	\end{theorem}

	\begin{proof}

		 In each of the cases in Theorems~\ref{thm:rect} and \ref{maxcross}, $A_r$ is a knot mosaic achieving the maximal predicted crossing number giving an upper bound for $c$.  First let's examine hexagonal 2 and 3-mosaics.  For hexagonal 2-mosaics we see the trefoil $\ddddot{A_2} = \hat{A_2} = \widehat{A_2}$ clearly achieves maximal crossing number for all three hexagonal settings.
		Moving up to $r=3$,  the hexagonal 3-mosaic cases $\ddddot{A_3}$, $\hat{A_3}$, and $\widehat{A_3}$ are all depicted in  Figure~\ref{fig:three}.
		 The exceptional case  is  standard hexagonal 3-mosaics.  We can smooth one crossing of $\ddddot{L_3}$ to get a knot, but it will always have a nugatory crossing so to get a reduced knot we must smooth 2 crossings showing $\ddddot{A_3}$ is maximal.  For semi-enhanced hexagonal 3-mosaics, we know they can only have crossings on interior tiles and thus can have crossing number at most 21.  It is easy to check that any 21 crossing semi-enhanced 3-mosaic is a link of at least 2 components so 
		 $\hat{A_3}$  with crossing number 20 must achieve the highest crossing number of any knot on a semi-enhanced 3-mosaic. The proof for $\widehat{A_3}$ is the same as the proof for $\widehat{A_r}$ $r>3$ so we will defer it to the paragraph below.

 Now we move on to the remaining cases: even rectangular boards with $r \geq 4$, standard and semi-enhanced hexagonal boards where  $r>3$, and enhanced hexagonal boards where $r > 2$.		 
		 In all of the remaining cases let $j$ be the number of components in $L_r$.  In each case $A_r$ is obtained by smoothing $j-1$ crossings of the saturated mosaic $L_r$.
		We know from Corollary~\ref{cor:smooth} that the knot which achieves the maximal crossing number is always achieved by smoothing crossings of a saturated mosaic.  In the even rectangular case there are two saturated mosaics up to crossing changes and $L_r$ requires fewer smoothings to become a reduced knot than the one with nugatory crossings in each of the four corners. In the hexagonal settings $L_r$ was always chosen to have the smallest number of components possible for any saturated link mosaic in its class, thus in all of these cases if we smooth fewer than $j-1$ crossings of a saturated mosaic we are left with a link of several components so it is impossible to find a knot with more crossings than $A_r$ on the mosaic board.  Since $A_r$ is reduced and alternating we know it achieves its crossing number in this embedding finishing the proof.

	\end{proof}

\section{Open Questions And Conjectures}

Hexagonal mosaics are rich with open questions.  Most theorems for rectangular mosaics translate into interesting questions in the hexagonal setting.  Here are a few specific ideas.

\bigskip

Ganzell and Henrich 
 show that any virtual knot or link can be represented as a virtual mosaic and go on to provide related computational results in \cite{gh}.  

\bigskip

\noindent Conjecture 1. \emph{
Virtual knots and links can be represented on virtual hexagonal mosaics.}

\bigskip

One would then hope to extend Ganzell and Henrich's computational results in this setting.

\bigskip

Alternatively, random knots are of great interest.  Mosaics are a good way to build random knots.  Hexagonal mosaics build knots far more efficiently than rectangular mosaics.  A rectangular 4-mosaic has at most 4 crossings and the only knots it supports are the unknot and the trefoil.  A hexagonal 4-mosaic can have up to 57 crossings and as shown above can contain knots of crossing number up to 55.

\bigskip

\noindent Conjecture 2. \emph{ (Easier Conjecture)
All knots can be created as a hexagonal $r$-mosaic  for some $r$ with interior tiles consisting only of the four tiles that contain 3 crossings and the two tiles that contain three arcs but no crossings.  
}

\bigskip

\noindent Conjecture 3. \emph{ (Harder Conjecture)
All knots can be created as a hexagonal $r$-mosaic  for some $r$ with interior tiles consisting only of the four 3 crossing tiles. }

\bigskip

Note that one uses five rectangular mosaic tiles (up to rotation) to generate every knot in rectangular mosaics, so this would imply that our computational complexity rises a little, but the number of knots we can create rises faster.

\bigskip

Finally, given a fixed $r$ this paper gave an upper bound for crossing number on a rectangular or hexagonal $r$-mosaic.  For rectangular $r$-mosaics  Lee, Hong, Lee, and Oh establish a bound in the opposite direction in \cite{Oh2b} showing that given a fixed crossing number, the mosaic number of a knot is less than or equal to that crossing number plus one.

\bigskip

\noindent Question 2. \emph{
	What is the best bound one can get when extending the results of \cite{Oh2b} to hexagonal $r$-mosaics?  }

\bigskip\bigskip


\bibliographystyle{amsplain}

  \bibliography{Hex}

\end{document}